\documentclass[a4paper]{amsart}
\usepackage[left=1.49in, right=1.49in, top=1.5in, bottom=1.8in]{geometry}
\newtheorem{theorem}{Theorem}[section]
\newtheorem{lemma}[theorem]{Lemma}

\theoremstyle{definition}

\theoremstyle{remark}

\numberwithin{equation}{section}

\newcommand{\abs}[1]{\lvert#1\rvert}

\usepackage[dvipsnames]{xcolor} 
\usepackage{hyperref}
\hypersetup{breaklinks, colorlinks,
	linkcolor = {Blue},
	citecolor = {Blue},
	urlcolor  = {OliveGreen}
}
\newcommand{\urlfa}[1]{\hypersetup{urlcolor = {Cyan%
	}}\url{#1}\hypersetup{urlcolor = {OliveGreen}}}%
\newcommand{\urlro}[1]{\hypersetup{urlcolor = {YellowOrange
	}}\url{#1}\hypersetup{urlcolor = {OliveGreen}}}%
\newcommand{\urlp}[1]{\hypersetup{urlcolor = {BrickRed%
	}}\url{#1}\hypersetup{urlcolor = {OliveGreen}}}%

\usepackage{scalerel}
\newcommand{\pr}[1]{\scalerel*{\textbf{(}}{\strut}#1\scalerel*{\textbf{)}}{\strut}}

\usepackage{enumitem}
\newlist{parts}{enumerate}{5}       
\setlist[parts]{wide, label=\textbf{\upshape(\alph*)}, ref=(\alph*), align=left, labelindent=0pt, labelsep=.5em, itemsep=.3ex, topsep=0ex, 
}

\let\oldsqrt\sqrt
\def\sqrt{\mathpalette\DHLhksqrt}
\def\DHLhksqrt#1#2{%
	\setbox0=\hbox{$#1\oldsqrt{#2\,}$}\dimen0=\ht0
	\advance\dimen0-0.2\ht0
	\setbox2=\hbox{\vrule height\ht0 depth -\dimen0}%
	{\box0\lower0.4pt\box2}}
	
\newcommand{\C}{{\mathbb C}} 

\newcommand{\N}{{\mathbb N}}

\newcommand{\R}{{\mathbb R}}

\newcommand{\cF}{{\mathcal F}}
\newcommand{\cG}{{\mathcal G}}

\newcommand{\cM}{{\mathcal M}}
\newcommand{\cP}{{\mathcal P}}

\newcommand{\rme}{\mathrm{e}}
\newcommand{\rmK}{\mathrm{K}}
\newcommand{\rmN}{\mathrm{N}}

\renewcommand{\epsilon}{\varepsilon}\renewcommand{\phi}{\varphi}
\renewcommand{\rho}{\varrho}
\renewcommand{\theta}{\vartheta}

\newcommand{\Eta}{\mathrm{H}}

\usepackage[outline]{contour}                                
\newcommand*{\fancy}[1]{{\color{white}\contour{black}{#1}}}  %

\newcommand{\norm}[1]{\left\rVert #1 \right\rVert}
\DeclareMathOperator*{\bigconv}{\mbox{\LARGE$\ast$}}
 
\newcommand{\dd}{{\mathrm{d}}}

\DeclareMathOperator*{\esssup}{ess\,sup}
\newcommand{\1}{\fancy{$1$}}  

\newcommand{\sgn}{{\mbox{\rm sgn}}}

\newcommand{\Prob}{\mbox{\rm Prob}} 

\DeclareFontFamily{U}{mathb}{\hyphenchar\font45}
\DeclareFontShape{U}{mathb}{m}{n}{
<-6> mathb5 <6-7> mathb6 <7-8> mathb7
<8-9> mathb8 <9-10> mathb9
<10-12> mathb10 <12-> mathb12
}{}
\DeclareSymbolFont{mathb}{U}{mathb}{m}{n}
\DeclareMathSymbol{\llcurly}{\mathrel}{mathb}{"CE}
\DeclareMathSymbol{\ggcurly}{\mathrel}{mathb}{"CF}

\usepackage{amsmath,mathtools,amssymb,amsthm,amscd,bm,graphicx}

\usepackage[T1]{fontenc}
\usepackage[english]{babel}

\begin{document}

\title[A sharper Lyapunov-Katz bound for summands Zolotarev-close to normal]{A sharper Lyapunov-Katz central limit error bound \\ for i.i.d.~summands Zolotarev-close to normal}

\author{Lena Jonas}
\address{Universit\"at Trier, Fachbereich IV -- Mathematik, 54286~Trier, Germany}
\email{jonasl@uni-trier.de}

\author{Lutz Mattner}
\address{Universit\"at Trier, Fachbereich IV -- Mathematik, 54286~Trier, Germany}
\email{mattner@uni-trier.de}

\subjclass[2020]{Primary 60F05, 60E15.}

\date{\today}


\keywords{Central limit theorem, sums of independent random variables}

\begin{abstract}
	We prove a central limit error bound for convolution powers of laws with finite moments of
	order~$r \in \mathopen]2,3\mathclose]$, taking a closeness of the laws to normality into account.
	Up to a universal constant,
	this generalises the case of~$r=3$ of the sharpening of the Berry~(1941) - Esseen~(1942) theorem
	obtained by~Mattner~(2024), namely by sharpening here the  Katz~(1963) error bound
	for the i.i.d.~case of Lyapunov's~(1901) theorem.
	
	Our proof uses a partial generalisation of the theorem of Senatov and Zolotarev (1986)
	used for the earlier special case.
	A result more general than our main one could be obtained by using
	instead another theorem of Senatov (1980),  but, unfortunately,
	an auxiliary inequality used in the latter's proof is wrong.
\end{abstract}

\maketitle
\section{Introduction and main results}            \label{sec:Intro}

The purpose of this paper is to prove
Theorem~\ref{Thm:Lyapunov-Katz_for_summands_Z-close_to_normal},
stated on page \pageref{page:Theorem_1.1} below,
which provides for the error in the central limit theorem,
for the convolution powers of a law with a finite moment of order
$r=2+\delta\in\mathopen]2,3\mathclose]$,
a bound taking a closeness of the initial law to normality into account.
Up to a universal constant, this generalises the special case of $r=3$
from Mattner~(2024)~\cite[p.~59, Theorem~1.5]{Mattner2024},
namely by improving the Katz~(1963)~\cite{Katz1963} error bound~\eqref{Eq:Lyapunov-Katz}
for the i.i.d.~case of Lyapunov's~(1901)~\cite{Lyapunov1901} theorem.
The improvement consists in replacing the $r$th standardised absolute moment
$\nu_{2+\delta}(\widetilde{P})$
occurring in~\eqref{Eq:Lyapunov-Katz} with a weak norm distance
of the initial law to normality,
namely $\Big( \zeta_{1} \vee \zeta_{2,\delta} \Big) \big( \widetilde{P} - \rmN \big)$
occurring in~\eqref{Eq:Lyapunov-Katz_for_summands_Z-close_to_normal}.
Here Senatov's~(1980)~\cite{Senatov1980} norm $\zeta_{2,\delta}$ as defined in
(\ref{Eq:Def_zeta_2,delta}, \ref{Eq:Def_zeta_m,g})
below satisfies $\zeta_{2,\delta}\le \zeta_{2+\delta}$,
where $\zeta_s$ with $s\in\mathopen]0,\infty\mathclose[$ denotes
the more common norm from~\eqref{Eq:Def_zeta_s}
introduced by Zolotarev~(1976) in~\cite[p.~725]{Zolotarev1976a} and~\cite[p.~386]{Zolotarev1976b},
for which Mattner and Shevtsova (2019)~\cite[pp.~497--498, 513--519]{MattnerShevtsova2019}
and \cite[pp.~56--58, 83--92]{Mattner2024} provide some references and elaborations,
and for which  $\zeta_{2,1} = \zeta_3$ in the earlier special case of $r=3$.

For stating our results and recalling some context more precisely,
let us first introduce some notation explained in
a bit more detail in~\cite[pp.~46--49]{Mattner2024}.
\pr{We also use, as in this sentence, parantheses for a few remarks or explanations
which could be skipped like footnotes.}
Besides the abbreviation {\em i.c.f.}~for {\em ignoring constant factors},
and the more common ones {\em i.i.d.}, {\em a.e.}, {\em w.l.o.g.},
{\em w.r.t.}, {\em iff},
$\text{\rm L.H.S.}$, $\text{\rm R.H.S.}$,
we use here also {\em i.a.d.}~for {\em independent but otherwise arbitrarily distributed},
to avoid the common but usually inaccurate ``non-i.i.d.''.
Order theoretic notions, like ``positive'' or ``increasing'', are understood in their wide sense.
We write $\cM$ for the algebra of all bounded signed measures on the Borel sets of $\R$,
with the inner multiplication being convolution,
$\big(M_1\ast M_2\big)(B)\coloneqq\iint \1_B(x+y)\,\dd M_1(x)\dd M_2(y)$
for $B\subseteq\R$ Borel. We accordingly write convolution powers as
$M^{\ast n} \coloneqq \bigconv_{j=1}^n M$
for $n\in\N_0\coloneqq\{0\}\cup\N$\,, with $M^0\coloneqq \delta_0$, the Dirac measure at $0\in\R$.
For $M\in\cM$, $|M|$ denotes its variation measure.
The subset of $\cM$ consisting of all probability measures, or \emph{laws},
is denoted by $\Prob(\R)$.
For $r\in[0,\infty[$ and $k\in\N_0$ we put
\begin{align*}
	\nu_r(M)  &\,\coloneqq\, \int|x|^r\,\dd|M|(x) \quad\text{for }M\in\cM\,,
	\quad \cM_r  \,\coloneqq\, \{ M\in\cM : \nu_r(M)<\infty\}\,,\\
	\mu_k(M)  &\,\coloneqq\, \int x^k\,\dd M(x) \quad\text{for }M\in\cM_k\,,
	\quad \mu \,\coloneqq\, \mu_1\,,  \\
	\cM_{r,k} &\,\coloneqq\, \{M\in\cM_r : \mu_0(M)=\ldots=\mu_k(M)=0\}
	\quad\text{if }k\le r\,, \\
	\sigma(P) &\,\coloneqq\, \sqrt{\textstyle{\iint \frac{1}{2}(x-y)^2\,\dd P(x)\dd P(y) } }
	\quad\text{for }P\in\Prob(\R)\,,\\
	\Prob_r(\R)&\,\coloneqq\, \Prob(\R)\cap\cM_r\,,
	\quad \cP_r\,\coloneqq\, \{P\in  \Prob_r(\R) : \sigma(P) >0\}\,.
\end{align*}
If $M\in\cM$, then $F_M(x)\coloneqq M(\,]-\infty,x]\,)$ for $x\in\R$ defines
its distribution function $F_M$, and if even $M\in\cM_{0,0}$, that is, $M\in\cM$
with $M(\R)=0$, then we call
$\left\|M\right\|_{\mathrm{K}}\coloneqq \sup_{x\in\R}\left|F_M(x)\right|$
the {\em Kolmogorov norm} of $M$.
If $P\in\cP_2$, then $\widetilde{P}$ denotes its standardisation, that is,
$\widetilde{P}(B)\coloneqq P(\sigma(P)B+\mu(P))$ for $B\subseteq \R$ Borel;
so $\widetilde{P}$ is the law of $(X-\mu(P))/\sigma(P)$
if $X$ is a random variable with the law $P$.
We accordingly put $\widetilde{\cP}\coloneqq \{\widetilde{P}:P\in\cP\}$ for $\cP\subseteq\cP_2$.
We write $\mathrm{N}_\sigma$
for the normal law with mean zero and standard deviation $\sigma$,
and $\mathrm{N}\coloneqq\mathrm{N}_1$ for the standard normal law.  
With this notation, the classical i.i.d.~central limit
theorem, proved in this generality first by Lindeberg~(1922)~\cite[p.~219, Satz~III]{Lindeberg1922},
can be written as
\begin{align}      \label{Eq:Lindeberg_CLT_iid}
	\lim_{n\rightarrow\infty}\big\|\widetilde{P^{\ast n}}-\mathrm{N}\big\|^{}_\mathrm{K}
	&\,=\, 0 \quad\text{for }P\in\cP_2\,,
\end{align}
and the present paper is one of the many aiming at bounding the approximation
error $\big\|\widetilde{P^{\ast n}}-\mathrm{N}\big\|^{}_\mathrm{K}$ for $n$ finite,
here in Theorem~\ref{Thm:Lyapunov-Katz_for_summands_Z-close_to_normal}
on page \pageref{page:Theorem_1.1}.

For a first reading of Theorem~\ref{Thm:Lyapunov-Katz_for_summands_Z-close_to_normal}, 
knowledge of the definitions~(\ref{Eq:Def_zeta_m,g}, \ref{Eq:Def_zeta_s}, \ref{Eq:Def_zeta_2,delta}) should suffice. 
The remaining text from here up to Theorem~\ref{Thm:Lyapunov-Katz_for_summands_Z-close_to_normal}
aims at providing some context, in particular the classical error bound~\eqref{Eq:Lyapunov-Katz},
and some understanding of the $\zeta$~norms used in the present paper.

For the purpose of bounding $\big\|\widetilde{P^{\ast n}}-\mathrm{N}\big\|^{}_\mathrm{K}\,$,
Katz~(1963)~\cite{Katz1963} was apparently the first to introduce
essentially as follows certain variants of the absolute moments $\nu_r$.
We put
\begin{align*}
	\nu_{m,g}(M)  &\coloneqq \int\!|x|^mg\big(|x|\big)\,\dd|M|(x)
	   \,\ \text{for }M\in\cM\,,
	\quad \cM_{m,g}  \coloneqq \{ M\in\cM : \nu_{m,g}(M)<\infty\}\,,
\end{align*}
for $m\in[0,\infty[$ and
$g:\mathopen[0,\infty\mathclose[ \rightarrow\mathopen[0,\infty\mathclose[$
measurable, but below there occurs only the case of $m\in\N_0$ and of $g$
belonging to the class $\cG$ of functions introduced as follows.

Let $\cG$ be the set of all functions $g:\mathopen[0,\infty\mathclose[ \rightarrow\mathopen[0,\infty\mathclose[$
which are increasing, with $\mathopen]0,\infty\mathclose[ \ni u\mapsto\frac{g(u)}{u}$ decreasing,
and which are not identically zero.
If $g\in\cG$ satisfies $g(1)=1$, 
then we call $g$ {\em normalised}. The set $\cG$ is a convex cone and a lattice,
with respect to the usual pointwise operations, so in particular for $g,g_1,g_2\in\cG$
also the {\em normalisation} $g/g(1)$ and the pointwise infimum $g_1\wedge g_2$ belong to $\cG$.
Examples of normalised functions $ g\in\cG$ occurring below are:
the identity $(\cdot)$, more generally
for $\delta\in\mathopen[0,1\mathclose]$
the power function $(\cdot)^\delta$ and its infimum  with the identity
$(\cdot)\wedge (\cdot)^\delta$,
and for $b\in[1,\infty[$ the function $g^{}_b$ given by
\begin{equation}              \label{Eq:Def_g_b}
	g^{}_b(u) \,\coloneqq\, u\wedge b \quad\text{ for }u\in[0,\infty[\,.
\end{equation}
Each $g\in\cG$ is subadditive,            \label{page:g_subadditive}
since for $u,v>0$ we have
$g(u+v)=(u+v)\frac{g(u+v)}{u+v}\le (u+v)\left(\frac{g(u)}{u} \wedge \frac{g(v)}{v}\right)
\le g(u)+g(v)$,
and obeys the scaling inequality $g(au)\le (1\vee a)g(u)$ for $a,u\in[0,\infty[\,$.

By the previous sentence, for $g\in\cG$,
the set $\cP_{2,g}\coloneqq \{ P\in\cP_2 : \nu_{2,g}(P)<\infty\}\subseteq\cM_{2,g}$
contains a $P\in\cP_2$ iff it contains its standardisation $\widetilde{P}$.
Katz~(1963)~\cite{Katz1963} proved,
imposing but not using the additional assumption of unboundedness of the functions $g$, that
\begin{align} \label{Eq:Katz_general_g}
	\norm{ \widetilde{P^{\ast n}} - \rmN }_{\rmK}
	\, \leq \,  \frac{c}{g(\sqrt{n})} \nu_{2,g}  \big( \widetilde{P} \big) \quad \text{for $g \in \cG$, $P \in \cP_{2,g}$, and $n \in \N$}
\end{align}
holds with some universal constant $c<\infty$,
with $c=1.8546$ being a possible choice
according to Korolev and Dorofeeva~(2017)
\cite[p.~46, Corollary~3, definition of $C_1$ on page 39]{KorolevDorofeeva2017}.
Since 
$\nu_{2,(\cdot)^\delta}=\nu_{2+\delta}$
for $\delta\in\mathopen[0,1\mathclose]$,
inequality \eqref{Eq:Katz_general_g}~obviously includes the result
\begin{align}                \label{Eq:Lyapunov-Katz}
	\norm{ \widetilde{P^{\ast n}} - \rmN }_{\rmK}
	\, \leq \,  \frac{c}{n^{\delta/2}_{}} \nu_{2+\delta}  \big( \widetilde{P} \big) \quad \text{for $\delta \in \mathopen[ 0, 1 \mathclose]$, $P \in \cP_{2+\delta}$, and $n \in \N$} \,,
\end{align}
with the same or perhaps some smaller universal constant $c$,
and~\eqref{Eq:Lyapunov-Katz} is what is improved
i.c.f.~to~\eqref{Eq:Lyapunov-Katz_for_summands_Z-close_to_normal}
in our main result,
Theorem~\ref{Thm:Lyapunov-Katz_for_summands_Z-close_to_normal} below.
In~\eqref{Eq:Lyapunov-Katz} we have included here the trivial case of $\delta=0$,
since its improvement in~\eqref{Eq:Lyapunov-Katz_for_summands_Z-close_to_normal}
is apparently not a trivial result.

How to name inequality~\eqref{Eq:Lyapunov-Katz}?
Restricted to $\delta\in\mathopen]0,1\mathclose[$ and
with unbounded constants $c_\delta<\infty$ in place of $c$,
it was, in a more general i.a.d.~version,
proved by Lyapunov~(1901)~\cite{Lyapunov1901} according to Petrov~(1995)~\cite[p.~184]{Petrov1995},
but hardly stated in a very obvious way in Lyapunov's paper.
To justify our claim of nonobviousness,
we refer to \cite[p.~18, $|\Delta|<\Omega$, $\Omega=\ldots$]{Lyapunov1901},
and observe that even Esseen~(1945)~\cite[p.~43, Theorem 1 restricted to $k<3$]{Esseen1945}
states a result seemingly weaker due to an additional second summand,
but actually equivalent as noted by Shevtsova~(2017)~\cite[p.~48]{Shevtsova2017}.
\pr{In more detail: If in Esseen's bound \cite[p.~43, (14)]{Esseen1945}
the first summand is $<$~$1$, then the second is even smaller
due to $\frac{1}{k-2}>1$. If the first summand is $\ge$~$1$,
then it is trivially larger than the left hand side of the inequality.
Hence the second summand may be omitted if Esseen's $c(k)$
is replaced with $\big(2\,c(k)\big)\vee 1$.}
Katz~(1963)~\cite{Katz1963} seems to be the first to prove with~\eqref{Eq:Katz_general_g}
a result obviously containing~\eqref{Eq:Lyapunov-Katz};
he reduces~\eqref{Eq:Katz_general_g} to the Berry-Esseen theorem, that is,
to the special case $\delta=1$ of~\eqref{Eq:Lyapunov-Katz}, acknowledging
help of his referee for obtaining a constant not depending on $g$,
yielding also a universal constant in~\eqref{Eq:Lyapunov-Katz}.
Hence one might name~\eqref{Eq:Lyapunov-Katz} after Lyapunov, Berry, Esseen, Katz,
and the latter's referee, but we have chosen here
``Lyapunov-Katz central limit error bound'' for short.

For more information on the history of the above results and some of its ramifications,
such as extensions to i.a.d.~situations,
better $\delta$-dependent constants in~\eqref{Eq:Lyapunov-Katz},
or analogues for so-called non-uniform error bounds, we recommend
the introductory pages of \cite{Shevtsova2017} as a starting point.

As Shevtsova~(2010) \cite[p.~226]{Shevtsova2010} points out,
the rate $n^{-\delta/2}$
in inequality~\eqref{Eq:Lyapunov-Katz}  is optimal
if just $\nu_{2+\delta}(\widetilde{P})$ is known:
Although, by 
Osipov and Petrov~(1967)
\cite[p.~283, Corollary]{OsipovPetrov1967}, we have
$\lim_{n\rightarrow\infty}n^{\delta/2}\text{L.H.S.\eqref{Eq:Lyapunov-Katz}}=0$
for $\delta \in\mathopen[0,1\mathopen[$ and every fixed $P\in\cP_{2+\delta}$,
in contrast to the case of $\delta=1$ where, by 
Esseen~(1956)~\cite[p.~162]{Esseen1956},
the limit in question also always exists but vanishes iff $P$ is a non-lattice law with
$\mu_3(\widetilde{P})=0$, we have
$\varliminf_{n\rightarrow\infty}n^{\delta/2}\sup_{P\in\cP_{2+\delta}}\text{L.H.S.\eqref{Eq:Lyapunov-Katz}}
 / \nu_{2+\delta}(\widetilde{P}) >0$ for every $\delta\in[0,1]$.
Here the final statement follows for  $\delta=0$
already from Matskyavichyus~(1983) \cite[p.~597, Theorem 2]{Matskyavichyus1983},
and follows for the remaining cases
from Shevtsova~(2010)~\cite{Shevtsova2010}.
\pr{In \cite[Theorem~3 on p.~243, definition of
$\overline{C}_\mathrm{AE}(\delta)$ on p.~238]{Shevtsova2010}
an explicit strictly positive lower bound is stated
with $\varlimsup$ instead of $\varliminf$, but proved with $\varliminf$ along even $n=2k$,
by taking the $n$th convolution power of a certain law $P_k$.
The result with  $\varliminf$ along arbitrary~$n$
then follows i.c.f.~by taking the $n$th convolution power of $P_n^{\ast 2}$.}

However, each of the upper bounds in~\eqref{Eq:Katz_general_g} and~\eqref{Eq:Lyapunov-Katz}
shares with the special Berry-Esseen case $\delta=1$ of~\eqref{Eq:Lyapunov-Katz}
the defect of  not becoming small for $n$ fixed but $P$ close to normal, since we have
\begin{align}                           \label{Eq:23/27}
	\nu_{2,g}(\widetilde{P}) &\,>\, \tfrac{23}{27}
	\quad\text{ for $g\in\cG$ normalised and $P\in\cP_2$\,},
\end{align}
as proved below on page~\pageref{Proof-23/27}, 
and even $\nu_{2+\delta}(\widetilde{P})\ge 1$ by, say, Jensen's inequality as
in \cite[p.~51, (1.27)]{Mattner2024}. Hence, ideally,
we would like to replace $\nu_{2,g}  \big( \widetilde{P} \big)$
in the upper bound in~\eqref{Eq:Katz_general_g}
with some weak norm distance to normality,
as achieved in \cite[Theorem~1.5]{Mattner2024} with
Zolotarev's norm distance $\big(\zeta_1\vee\zeta_3\big)(\widetilde{P}-\mathrm{N})$
for the special Berry-Esseen case of $\delta=1$ in~\eqref{Eq:Lyapunov-Katz},
for $n\ge 2$.
To this end, we consider here norm distances more general than \mbox{$\zeta_1$ and $\zeta_3$},
essentially following  Senatov~(1980)~\cite{Senatov1980}:

For $m \in \N_0$ and $g\in\cG$, let $\cF_{m,g}$ be the set of all $m$ times
differentiable measurable functions $f:\R\rightarrow \R$ for which
their $m$th derivatives $f^{(m)}$ satisfy the increment boundings
\begin{align}               \label{Eq:Increment_bound_for_mth_derivative_with_g}
	\left| f^{(m)}(x)- f^{(m)}(y) \right|
	&\,\le\, g\big(\,|x-y|\,\big) \quad\text{ for }x,y\in\R\,,
\end{align}
and let $\cF_{m,g}^\infty$ be the set of all bounded $f\in\cF_{m,g}$.
Then
\begin{equation}                          \label{Eq:Def_zeta_m,g}
	\zeta_{m,g}(M) \,\coloneqq\, \sup\left\{
	\,{\textstyle \left| \int f\,\dd M\right|} : f \in \cF^\infty_{m,g} \,\right\}
	\quad \text{ for }M\in\cM
\end{equation}
defines on the vector space $\cM$ an {\em enorm}
(by which we mean a norm except that the value $\infty$ is also allowed),
which, since $\cF_{m,g}^\infty$ is a translation invariant
set of bounded measurable functions,  satisfies the so-called
regularity inequality
\begin{align}               \label{Eq:regularity}
	\zeta_{m,g}(M_1\ast M_2) &\,\le\, \zeta_{m,g}(M_1)\nu_0(M_2)
	\quad\text{ for }M_1,M_2\in\cM
\end{align}
and its variants listed in
\cite[p.~83, Lemma~5.1 up to the line after (5.7)]{Mattner2024}.
We have the obvious implications
\begin{align}              \label{Eq:zeta_g_1-g_2_comparison_inequality}
	m\in\N_0\,,\ g^{}_1,g^{}_2\in\cG\,,\ g^{}_1\le g^{}_2
	 &\ \Rightarrow\ \zeta_{m,g^{}_1} \le  \zeta_{m,g^{}_2}\,, \\
	m\in\N_0\,,\ g\in\cG\,,\     \label{Eq:zeta_scaling_w.r.t._g}
	    \lambda\in\mathopen]0,\infty\mathclose[
	 &   \ \Rightarrow\ \zeta_{m,\lambda g} = \lambda \zeta_{m, g}\,.
\end{align}
The most common special cases of the $\zeta_{m,g}$
from~\eqref{Eq:Def_zeta_m,g} are given by
\begin{align}      \label{Eq:Def_zeta_s}
	\zeta_{m+\delta} &\,\coloneqq\, \zeta_{m,(\cdot)^\delta}
	\quad\text{ for $m\in\N_0$ and $\delta\in\mathopen]0,1\mathclose]\,$},
\end{align}
which satisfy, in addition to~\eqref{Eq:regularity}, also a scaling identity,
often called {\em homogeneity},
\begin{equation}                 \label{Eq:scaling_identity_zeta_s}
 \zeta_s(M(\tfrac{\cdot}{a})) \,=\, a^s\zeta_s(M)
   \quad\text{ for }s,a\in\mathopen]0,\infty\mathclose[
   \text{ and } M\in\cM\,,
\end{equation}
as a special case of \cite[pp.~83--84, (5.8) implies (5.9)]{Mattner2024}.

The definition of $\zeta_s$ for $s\in\mathopen]0,\infty\mathclose[$
given by~\eqref{Eq:Def_zeta_s} conforms with the one Zolotarev~\cite{Zolotarev1976a,Zolotarev1976b}
gave for the special case  of signed measures
being differences of laws, $M=P-Q$ with $P,Q\in\Prob(\R)$,
writing then $\zeta_s(X,Y)$ with random variables $X, Y$ with laws $P,Q$,
rather than $\zeta_s(P-Q)$ as here and in~\cite{Mattner2024},
or $\zeta_s(P,Q)$ as in~\cite{MattnerShevtsova2019} and in several other publications.
The particular case of $\zeta_1(P-Q)$
is the Kantorovich distance between $P$ and~$Q$, often rather called
Wasserstein distance, for which \cite[pp.~54--55]{Mattner2024}
gives some historical notes.
For normalised but otherwise arbitrary functions $g\in\cG$,
Senatov~(1980)~\cite[p.~746]{Senatov1980} introduced as~$\zeta_g(P,Q)$
what in our present notation is $\zeta_{2,(\cdot)\wedge g}(P-Q)$,
except that
he forgot to introduce a boundedness condition on the functions $f$
ensuring well-definedness of $\int f\,\dd|M|$
in~\eqref{Eq:Def_zeta_m,g} with $m=2$ and $(\cdot)\wedge g$ in place of $g$.
The definition of $\zeta_{m,g}(M)$ in~\eqref{Eq:Def_zeta_m,g} remains unaffected
by assuming, as for example in~\cite{Mattner2024}, the functions in $\cF_{m,g}^\infty$
to be $\C$-valued instead of $\R$-valued,
since for any such function $f$ then
$\int f\,\dd M = cr $ with $c\in\C,|c|=1, r\in[0,\infty[\,$,
and hence $\left|\int f\,\dd M\right| = r = \int c^{-1}f\,\dd M = \int \Re( c^{-1}f)\,\dd M$
with the real part $\Re( c^{-1}f)$ belonging to $\cF_{m,g}^\infty$.
With this remark, and using the fact that a function $f:\R\rightarrow\R$
satisfies $|f(x)-f(y)|\le |x-y|$ for $x,y\in\R$ iff it is absolutely
continuous with $\esssup_{x\in\R}|f'(x)|\le 1$, one easily checks that
for $r\in\N$ the definitions of~$\zeta_r$ here and in
\cite[p.~57, Definition 1.3]{Mattner2024} are indeed equivalent.

In contrast to~\eqref{Eq:zeta_g_1-g_2_comparison_inequality},
for example $\zeta_1$ and $\zeta_3$ are incomparable
even at arguments being of special interest in the present context,
namely $M=\widetilde{P}-\mathrm{N}$ close to the zero measure
and $P\in\widetilde{\cP_2}$ quite natural, as in
\cite[Examples 1.7 and 1.9]{Mattner2024}.

The main result of the present paper,
Theorem~\ref{Thm:Lyapunov-Katz_for_summands_Z-close_to_normal} just below,
and also Theorem~\ref{Thm--Z/S_zeta1_zeta2,delta} used for proving the former,
involve besides $\zeta_1$ the special cases of $\zeta_{2,g}$ given by
\begin{align}
	\zeta_{2,\delta} &\,\coloneqq\,                  \label{Eq:Def_zeta_2,delta}
	\zeta_{2,(\cdot)\wedge(\cdot)^\delta}
	\quad\text{ for $\delta\in\mathopen[0,1\mathclose]$}\,,
\end{align}
which are by~\eqref{Eq:zeta_g_1-g_2_comparison_inequality} with $m=2$ increasing in their parameter,
\begin{align}                    \label{Eq:zeta_2,delta_ordered}
  0\le \delta\le\epsilon\le 1 &\ \Rightarrow\ \zeta_{2,\delta} \,\le\, \zeta_{2,\epsilon}\,,
\end{align}
and which can be upper-bounded at least on $\cM_{2,2}$
by the more classical corresponding enorms $\zeta_{2+\delta}$ through
\begin{align}                \label{Eq:zeta_2,delta_vs_zeta_2+delta}
	\zeta_{2,\delta} &\,\le\, \zeta_{2,(\cdot)^\delta}
	\,
	\left\{\begin{array}{ll}
		\!=\,\zeta^{}_{2+\delta} \phantom{\int\limits_0}
		&\text{ if }\delta\in\mathopen]0,1\mathclose]\,,\\
		\!\eqqcolon\, \zeta_2^\flat \,\le\, \zeta_2 &\text{ if }\delta = 0
	\end{array}\right\},
	\quad \zeta_{2,1}\,=\,\zeta_3\,,
\end{align}
with $\zeta_2^\flat \,\le\, \zeta_2$ on $\cM_{2,2}$
proved on page~\pageref{page:prof_of_zeta_2^flat_le_zeta},
and the other claims actually holding on all of~$\cM$,
by~\eqref{Eq:zeta_g_1-g_2_comparison_inequality} or just by definition.
In particular we have $\zeta_{2,\delta} \le\,\zeta_3\,$.

\begin{theorem}[Lyapunov-Katz for i.i.d.~summands Zolotarev-close to normal]
                        	\label{Thm:Lyapunov-Katz_for_summands_Z-close_to_normal}
                        	\label{page:Theorem_1.1}
	There exists a constant $c \in \mathopen] 0, \infty \mathclose[$ satisfying
	\begin{align}                \label{Eq:Lyapunov-Katz_for_summands_Z-close_to_normal}
		\norm{ \widetilde{P^{\ast n}} - \rmN }_{\rmK} \, \leq \, \frac{c}{n_{}^{\delta/2}}
		\Big( \zeta_{1} \vee \zeta_{2,\delta} \Big) \big( \widetilde{P} - \rmN \big) \quad \text{for $\delta \in \mathopen[ 0, 1 \mathclose]$, $P \in \cP_{2+\delta}$, and $n \geq 2$} \,.
	\end{align}
	One may take here $c = 48$.
\end{theorem}

I.c.f., Theorem~\ref{Thm:Lyapunov-Katz_for_summands_Z-close_to_normal} generalises
the earlier special case \cite[p.~59, Theorem 1.5]{Mattner2024} of \mbox{$\delta=1$},
and sharpens the Lyapunov-Katz theorem~\eqref{Eq:Lyapunov-Katz} as desired,
since for $\delta\in[0,1]$ and \mbox{$P\in\cP_{2+\delta}$} we may
use~\eqref{Eq:zeta_2,delta_vs_zeta_2+delta},
Lemma~\ref{Lem:cF_m,g-bounds_and_zeta_m,g_for_unbounded_f}\ref{part:zeta_s_le_nu_s},
that $\log(\nu_s(Q))$ is convex in $s\in[0,\infty[$ and vanishes on $\{0,2\}$ if $Q\in\widetilde{\cP}$,
$\nu_1(\rmN)=\frac{2}{\sqrt{2\pi}}$, and
$\nu_{2+\delta}(\rmN)\le \nu_3(\rmN)=\frac{4}{\sqrt{2\pi}}$, to get
\begin{align}            \label{Eq:zeta_1-distance_to_normality_bounded}
 \zeta_1\big(\widetilde{P} - \rmN \big ) \,\le\, \nu_1\big(\widetilde{P} - \rmN \big)
  \,\le\,  \nu_1\big(\widetilde{P}\big) + \tfrac{2}{\sqrt{2\pi}}
  \,\le\,1+\tfrac{2}{\sqrt{2\pi}}
  \,\le\,\left(1+\tfrac{2}{\sqrt{2\pi}}\right)\nu_{2+\delta}\big(\widetilde{P}\big)
\end{align}
and
\begin{align}           \label{Eq:zeta_2,delta-distance_to_normality_bounded}
 \zeta_{2,\delta}\big(\widetilde{P} - \rmN \big )
   \,\le\, \zeta_{2+\delta}\big(\widetilde{P} - \rmN \big )
   \,\le\, \frac{ \nu_{2+\delta}\big(\widetilde{P} - \rmN \big )}{(2+\delta)(1+\delta)}
   \,\le\, \left(\tfrac{1}{2}+\tfrac{2}{\sqrt{2\pi}}\right) \nu_{2+\delta}\big(\widetilde{P}\big )\,,
\end{align}
hence
$
 \Big( \zeta_{1} \vee \zeta_{2,\delta} \Big) \big( \widetilde{P} - \rmN \big)
  \,\le\,\left(1+\tfrac{2}{\sqrt{2\pi}}\right)\nu_{2+\delta}\big(\widetilde{P}\big)
$. Here the final bounds for the left hand sides
in (\ref{Eq:zeta_1-distance_to_normality_bounded},~\ref{Eq:zeta_2,delta-distance_to_normality_bounded})
can sometimes be improved,
namely we have $\zeta_1\big(\widetilde{P} - \rmN \big ) \,\le\,\nu_{3}\big(\widetilde{P}\big)$
if $\delta=1$ by Goldstein (2010) \cite[Theorem~1.1]{Goldstein2010}
and by Tyurin (2010) \cite[Theorem~4]{Tyurin2010},
and
$\zeta_{2+\delta}\big(\widetilde{P} - \rmN \big )
\,\le\, \frac{1}{(2+\delta)(1+\delta)}\nu_{2+\delta}\big(\widetilde{P}\big)
\,<\, \frac{1}{2}\nu_{2+\delta}\big(\widetilde{P}\big)$
if $\delta>0$ by Tyurin~(2012) \cite[Theorem~2]{Tyurin2012},
with $\zeta_{3}\big(\widetilde{P} - \rmN \big )
\,\le\, \frac{1}{6}\nu_3(\widetilde{P})$ already by \cite[Theorem~4]{Tyurin2010}.

Again i.c.f., Theorem~\ref{Thm:Lyapunov-Katz_for_summands_Z-close_to_normal}
improves the following earlier sharpening of the
Lyapunov-Katz theorem~\eqref{Eq:Lyapunov-Katz} due to
Senatov~(1980)~\cite[p.~747, Theorem 1 for $k=1$, $g=(\cdot)^\delta$]{Senatov1980}:
There exists a constant $c\in\mathopen]0,\infty\mathclose[$ with
\begin{align} \label{Eq:Lyapunov-Katz-Senatov_for_summands_Z-close_to_normal}
	\norm{ \widetilde{P^{\ast n}} - \rmN }_{\rmK}  \, \leq \,  \frac{c}{n^{\delta/2}_{}}
	\Big(\norm{\,\cdot\,}_{\rmK} \vee  \zeta_{1} \vee \zeta_{2,\delta} \Big) \big( \widetilde{P} - \rmN \big) \quad \text{for $\delta \in \mathopen[ 0, 1 \mathclose]$, $P \in \cP_{2+\delta}$,  and
	$n \in \N$} \,.
\end{align}
Due to the additional Kolmogorov norm $\norm{\,\cdot\,}_\rmK$ on the right,
inequality~\eqref{Eq:Lyapunov-Katz-Senatov_for_summands_Z-close_to_normal}
is trivial for~$c\ge 1$ and $n=1$; hence~\eqref{Eq:Lyapunov-Katz_for_summands_Z-close_to_normal}
indeed improves~\eqref{Eq:Lyapunov-Katz-Senatov_for_summands_Z-close_to_normal}.

By~\eqref{Eq:zeta_2,delta_vs_zeta_2+delta} we may replace $\zeta_{2,\delta}$
in~\eqref{Eq:Lyapunov-Katz_for_summands_Z-close_to_normal} with $\zeta_{2+\delta}$.
But it appears likely to us that, for each fixed $\delta\in\mathopen]0,1[\,$,
one can have
$ \Big( \zeta_{1} \vee \zeta_{2+\delta} \Big) \big( \widetilde{P_\epsilon} - \rmN \big)
\ggcurly
\Big( \zeta_{1} \vee \zeta_{2,\delta} \Big) \big( \widetilde{P_\epsilon} - \rmN \big) \rightarrow 0$
for appropriate $P_\epsilon\in\cP_{2+\delta}$, similarly as in~\cite[p.~64, Example~1.12]{Mattner2024}.

The constant $c$ in Theorem~\ref{Thm:Lyapunov-Katz_for_summands_Z-close_to_normal}
can not be chosen smaller than $\sqrt{2}\frac{15+6\sqrt{3}}{13\sqrt{2\pi}}=1.1020\ldots$ $\eqqcolon$ $c_1$,
and for every fixed value of $\delta$ not smaller than $\sqrt{2}^{\delta-1} c_1$,
by the example \cite[p.~108, (12.16)]{Mattner2024}, since there
$\zeta_3$ is negligeable compared to $\zeta_1$, and hence so is $\zeta_{2,\delta}\le\zeta_3$ by~\eqref{Eq:zeta_2,delta_ordered} and~\eqref{Eq:zeta_2,delta_vs_zeta_2+delta}.

In the special case of $\delta=1$, \cite[p.~59, Theorem 1.5]{Mattner2024}
gives~\eqref{Eq:Lyapunov-Katz_for_summands_Z-close_to_normal} with $c=7.2\,$.

Theorem~\ref{Thm:Lyapunov-Katz_for_summands_Z-close_to_normal}, with a worse constant,
was first obtained by Jonas~(2024)~\cite{Jonas2024}, 
by combining~\eqref{Eq:Lyapunov-Katz-Senatov_for_summands_Z-close_to_normal}
with the convolution inequality \cite[p.~71, Corollary~3.3]{Mattner2024}.
\pr{In more detail, assuming for simplicity $n=2k$ even:
One applies~\eqref{Eq:Lyapunov-Katz-Senatov_for_summands_Z-close_to_normal}
to $(P^{\ast 2},k)$ rather than to $(P,n)$, then uses on the right
from the just cited corollary $\norm{ \widetilde{P^{\ast 2}} - \rmN }_{\rmK}
\le \frac{2}{\sqrt{2\pi}}\zeta_1(\widetilde{P} - \rmN )$,
and for the other two norms the regularity~\eqref{Eq:regularity}
applied to $M_1\coloneqq \widetilde{P}-\mathrm{N}$
and $M_2\coloneqq \widetilde{P}+\mathrm{N}$,
prepared with the scaling identity~\eqref{Eq:scaling_identity_zeta_s} for $\zeta_1$,
and with the scaling inequality \eqref{Eq:Scaling_ineq_for_zeta_2,delta}
for~$\zeta_{2,\delta}$.}

In the present paper, we instead first provide the
Senatov-Zolotarev type $\zeta_1\vee\zeta_{2,\delta}$
Theorem~\ref{Thm--Z/S_zeta1_zeta2,delta},
partially generalising the more complete result for $\delta=1$
obtained in \cite[p.~70, Theorem 3.1]{Mattner2024},
and there named after Zolotarev only, erroneously as explained in
section~\ref{sec:Errata} below.
Then we prove Theorem~\ref{Thm:Lyapunov-Katz_for_summands_Z-close_to_normal}
by combining Theorem~\ref{Thm--Z/S_zeta1_zeta2,delta}
with the convolution inequality from \cite{Mattner2024},
and by using the Lyapunov-Katz inequality~\eqref{Eq:Lyapunov-Katz}
in the case of large values of
$\zeta_{2,\delta}(\widetilde{P}-\mathrm{N})$.
This approach appears to be simpler overall
and leads to a better constant in
Theorem~\ref{Thm:Lyapunov-Katz_for_summands_Z-close_to_normal}.
Also, Theorem~\ref{Thm--Z/S_zeta1_zeta2,delta}, in spite of its imperfectness
apparent if compared with \cite[Theorem 3.1]{Mattner2024},
seems to be of some independent interest to us,
while~\eqref{Eq:Lyapunov-Katz-Senatov_for_summands_Z-close_to_normal}
is superseded by Theorem~\ref{Thm:Lyapunov-Katz_for_summands_Z-close_to_normal}.

\begin{theorem}[A Senatov-Zolotarev type $\zeta_{1} \vee \zeta_{2,\delta}$ Theorem]
	\mbox{}                                      \label{Thm--Z/S_zeta1_zeta2,delta}
	\begin{parts}
		\item \label{Thm--Z/S_zeta1_zeta2,delta_1stPart}
		Let $\delta \in \mathopen[ 0, 1 \mathclose]$. We have
		\begin{align} \label{Eq--Z/S_zeta1_zeta2,delta_1stPart}
			\zeta_{1}( \widetilde{P^{\ast n}} - \rmN )  \, \leq \,  \frac{1}{n^{\delta/2}_{}} \xi^{}_{\delta} \Big( \zeta_{1}( \widetilde{P} - \rmN ), \zeta_{2,\delta}( \widetilde{P} - \rmN ) \Big) \quad \text{for $P \in \cP_{2+\delta}$ and $n \in \N$} \,,
		\end{align}
		where the function $\xi^{}_{\delta} : \mathopen[0, \infty \mathclose[^{2} \to \mathopen[0, \infty \mathclose]$,
		with $\xi^{}_{\delta}(\varkappa,\zeta)=\infty$
		only if $\delta =0$ and $\zeta \ge (2\gamma_0)^{-1}=0.2582\ldots$,
		is defined through
		\begin{gather}
			\alpha^{}_{\delta}   \coloneqq   \left( \tfrac{16\rme^{-1/2} - 4}{\sqrt{2\pi}} \right)^{1-\delta} \left( \tfrac{4\mathrm{e}^{-1/2}}{\sqrt{2\pi}} \right)^{\delta},
			\quad \beta   \coloneqq   \tfrac{4}{\sqrt{2\pi}} \,, 
			\, \gamma^{}_{\delta}   \coloneqq   \left( \tfrac{8\rme^{-1/2}}{\sqrt{2\pi}} \right)^{1-\delta} \left( \tfrac{2 + 8\rme^{-3/2}}{\sqrt{2\pi}} \right)^{\delta} , \label{Eq--Z/S_zeta1_zeta2,delta_1stPart_alpha-beta-gamma} \\
			g^{}_{\delta}(\eta)  \, \coloneqq \,  2\frac{(1 + \eta^2)^{\frac{1-\delta}{2}}}{\eta} \quad \text{for $\eta \in \mathopen] 0, \infty \mathclose[$} \,,
			\label{Eq--Z/S_zeta1_zeta2,delta_1stPart_g}\\
			\xi^{}_{\delta}( \varkappa, \zeta )  \, \coloneqq \,  \inf\left\{ \frac{\varkappa + \alpha^{}_{\delta}\zeta + \beta\eta}{1 - \gamma^{}_{\delta} g^{}_{\delta}(\eta) \zeta} \, : \, \eta \in \mathopen] 0, \infty \mathclose[ \,, \, \gamma^{}_{\delta} g^{}_{\delta}(\eta) \zeta < 1 \right\} \;\:
			\text{for $\varkappa, \zeta \in \mathopen[0, \infty \mathclose[$} \,.
			\label{Eq--Z/S_zeta1_zeta2,delta_1stPart_xi}
		\end{gather}

		\item \label{Thm--Z/S_zeta1_zeta2,delta_2ndPart}
		Putting
		\begin{align}
			\alpha \coloneqq \frac{16\rme^{-1/2} - 4}{\sqrt{2\pi}}, 
			\quad \beta \coloneqq \frac{4}{\sqrt{2\pi}},
			\quad \gamma \coloneqq \frac{8\rme^{-1/2}}{\sqrt{2\pi}}
		\end{align}
		we have, for any choice of $(\theta, c) \in \mathopen] 0, 1 \mathclose[ \times \mathopen] 0, \infty \mathclose[$ satisfying
		$\frac{\theta^2}{4\gamma^2} - \frac{1}{c^2} > 0$, the following: For $\delta \in \mathopen[ 0, 1\mathclose]$ and $P$ satisfying $\zeta_{2,\delta}( \widetilde{P} - \rmN )
		\leq \sqrt{\frac{\theta^2}{4\gamma^2} - \frac{1}{c^2}}$ we have
		\begin{align} \label{Eq--Z/S_zeta1_zeta2,delta_2ndPart}
			\zeta_{1}( \widetilde{P^{\ast n}} - \rmN )  \, \leq \,  \frac{\, 1 + \alpha + \beta c \ }{1 - \theta} \frac{\Big( \zeta_{1} \vee \zeta_{2,\delta} \Big)\big( \widetilde{P} - \rmN \big)}{n^{\delta/2}_{}} \,.
		\end{align}
	\end{parts}
\end{theorem}

We finally have to mention that
Senatov~(1980)~\cite[p.~747, Theorem 1 for $k=1$]{Senatov1980}
states a bound more general than
\eqref{Eq:Lyapunov-Katz-Senatov_for_summands_Z-close_to_normal},
namely one in our present notation equivalent to
the existence of a constant $c\in\mathopen]0,\infty\mathclose[$ with
\begin{align} \label{Eq:Katz-Senatov_UNPROVEN_for_summands_Z-close_to_normal}
	\norm{ \widetilde{P^{\ast n}} - \rmN }_{\rmK}
	&\,\le\, \frac{c}{g(\sqrt{n})}
	\Big(\norm{\,\cdot\,}_{\rmK} \vee  \zeta_{1} \vee \zeta_{2,g} \Big) \big( \widetilde{P} - \rmN \big) \quad \text{for $g \in \cG\,$, $P \in \cP_{2,g}\,$, and $n \in \N$} \,,
\end{align}
which however seems to be neither proven nor disproven so far.

\pr{Senatov has, in our notation,
\eqref{Eq:Katz-Senatov_UNPROVEN_for_summands_Z-close_to_normal}
with $g$ restricted to be normalised and to satisfy $g=(\cdot)\wedge g$.
But if a normalised $g\in\cG$ is given, then Senatov's statement
applied to $(\cdot)\wedge g$ yields the inequality
in~\eqref{Eq:Katz-Senatov_UNPROVEN_for_summands_Z-close_to_normal}
for the given $g$, since
$\sqrt{n}\wedge g(\sqrt{n})=g(\sqrt{n})$ due to
$g(\sqrt{n})=\sqrt{n}\frac{g(\sqrt{n})}{\sqrt{n}}\le \sqrt{n}\frac{g(1)}{1}=\sqrt{n}$,
and since $\zeta_{2,(\cdot)\wedge g} \le \zeta_{2,g}$
by~\eqref{Eq:zeta_g_1-g_2_comparison_inequality}. And if finally $g\in\cG$ is arbitrary,
then the inequality
in~\eqref{Eq:Katz-Senatov_UNPROVEN_for_summands_Z-close_to_normal}
applied to the normalised function $\frac{g}{g(1)}$ yields the one for $g$,
due to~\eqref{Eq:zeta_scaling_w.r.t._g}.}

In his attempt at proving~\eqref{Eq:Katz-Senatov_UNPROVEN_for_summands_Z-close_to_normal},
Senatov uses a scaling inequality more general than
Lemma~\ref{Lem--Scaling_ineq_for_zeta2,delta} below,
which he states
in~\cite[p.~746, third display from below]{Senatov1980}
without any proof or reference, and which is refuted
by our Lemma~\ref{Lem:Senatov_scaling_inequality_wrong}.

If~\eqref{Eq:Katz-Senatov_UNPROVEN_for_summands_Z-close_to_normal} were true,
then we should get, for some constant $c\in\mathopen]0,\infty\mathclose[\,$,
\begin{align} \label{Eq:general_Katz_UNPROVEN_for_summands_Z-close_to_normal}
	\norm{ \widetilde{P^{\ast n}} - \rmN }_{\rmK}
	&\,\le\, \frac{c}{g(\sqrt{n})}
	\Big(\zeta_{1} \vee \zeta_{2,g} \Big) \big( \widetilde{P} - \rmN \big) \quad \text{for $g \in \cG\,$, $P \in \cP_{2,g}\,$, and
	$n\ge2$}\,,
\end{align}
similarly as Theorem~\ref{Thm:Lyapunov-Katz_for_summands_Z-close_to_normal}
was derived from~\eqref{Eq:Lyapunov-Katz-Senatov_for_summands_Z-close_to_normal}
in \cite{Jonas2024}. Applied to the functions $g_b$ from~\eqref{Eq:Def_g_b},
this would sharpen the rowwise i.i.d.~case of Lindeberg's theorem.

\section{Proofs of the two main results, Theorems~\ref{Thm:Lyapunov-Katz_for_summands_Z-close_to_normal} and~\ref{Thm--Z/S_zeta1_zeta2,delta}}
In the rest of this paper, we often abbreviate $M_1M_2\coloneqq M_1\ast M_2$,
and accordingly $M^n\coloneqq M^{\ast n}$ for $n\in\N_0$\,.
	
\begin{proof}[Proof of Theorem~\ref{Thm--Z/S_zeta1_zeta2,delta}]
\ref{Thm--Z/S_zeta1_zeta2,delta_1stPart} We follow
	\cite[pp. 93--94, proof of Theorem~3.1]{Mattner2024} in using the commutative ring identity
	obtained by the calculation, valid for $n\in\N$
	and $\sum_{j\in\emptyset}\ldots\coloneqq 0$,
	\begin{align} \label{Eq--ring_identity}
		P^{n} - Q^{n}  &\, =\,  \sum_{j=0}^{n-1} \left( P^{n-j}Q^{j} - P^{n-j-1}Q^{j+1} \right) \\
		&\, = \,  P^{n} - P^{n-1}Q \, + \, \sum_{j=1}^{n-1} \left( P^{n-j}Q^{j} - P^{n-j-1}Q^{j+1} \right) \nonumber\\
		&\, = \,  (P - Q)P^{n-1} \nonumber\\ &\qquad + \sum_{j=1}^{n-1} \Big( \left(PQ^{n-1} - Q^{n}\right) + \left(P^{n-j-1} - Q^{n-j-1}\right)\left(PQ^{j} - Q^{j+1}\right) \Big) \nonumber\\
		&\, = \,  (P - Q)P^{n-1} + (n-1)  (P-Q)Q^{n-1} \nonumber\\
		&\qquad  + \sum_{j=1}^{n-2} \Big( \left(P^{n-j-1} - Q^{n-j-1}\right)  (P-Q)Q^{j}  \Big) \,. \nonumber
	\end{align}

 Let now $\delta \in \mathopen[ 0, 1 \mathclose]$
 and $P \in \cP_{2+\delta}\,$,
 and let $\xi^{}_{\delta, 0} \coloneqq \xi^{}_{\delta} \big( \zeta_{1}( \widetilde{P} - \rmN ), \zeta_{2,\delta}( \widetilde{P} - \rmN ) \big)$. Let $n \in \N$ be such that we have
	\begin{align} \label{Eq--Z/S_inductive_hypothesis}
		\zeta_{1}( \widetilde{P^{\ast k}} - \rmN )
		 \, \leq \,  \frac{\xi^{}_{\delta, 0}}{k^{\delta/2}_{}} \quad \text{for $k \in \{1, \dots, n-2\}$} \,.
	\end{align}
	Using this inductive hypothesis, which, as in \cite[pp. 93--94]{Mattner2024}, is trivial for $n=1$ and does not involve $k=n-1$, we are going to prove
	\begin{equation} \label{Eq--Z/S_inductive_goal}
	   \zeta_{1}( \widetilde{P^{\ast n}} - \rmN ) \, \leq \,  \frac{\xi^{}_{\delta, 0}}{n^{\delta/2}_{}}\,.
	\end{equation}
	To this end, we may assume w.l.o.g. $P=\widetilde{P}$.
	With $Q \coloneqq \rmN$ we then have
	\begin{align}                        \label{Eq:zeta_1_scaled_proof_of_Th_1.2}
	 \zeta_1(\widetilde{P^{\ast n}} - \rmN)
	    \, = \, \frac{1}{\sqrt{n}} \zeta_1(P^{n} - Q^{n})
	\end{align}
	by the scaling identity~\eqref{Eq:scaling_identity_zeta_s},
	and we are going to apply~\eqref{Eq--ring_identity}
	to R.H.S.~\eqref{Eq:zeta_1_scaled_proof_of_Th_1.2}
	after employing a so-called smoothing-inequality,
	where $P^{n} - Q^{n}$ is convolved with 	the normal law~$\mathrm{N}_\epsilon$
	with mean zero and yet unspecified standard deviation $\epsilon$.

	For $\epsilon \in \mathopen] 0, \infty \mathclose[$ we thus get,
	recalling the notation~\eqref{Eq--Z/S_zeta1_zeta2,delta_1stPart_alpha-beta-gamma},
	\begin{align*}
		\sqrt{n}\,\zeta_{1}( \widetilde{P^{\ast n}} - \rmN )  &\, = \,  \zeta_{1}( P^{n} - Q^{n} ) \\
		&\,\leq \,  \beta\epsilon + \zeta_{1}\big( (P^{n} - Q^{n})\rmN_{\epsilon} \big) \\
		&\,\leq \,  \beta\epsilon + \zeta_{1}\big( (P - Q)\rmN_{\epsilon}P^{n-1} \big) + (n-1) \zeta_{1}\left( (P - Q)Q^{n-1}\rmN_{\epsilon} \right) \\ &\qquad + \sum_{j=1}^{n-2} \zeta_{1}\left( ( P^{n-j-1} - Q^{n-j-1} )  (P -Q )Q^{j}\rmN_{\epsilon}  \right) \\
		&  
		\, \leq \,  \beta\epsilon + \zeta_{1}( P - Q) + (n-1) \zeta_{1}\left( (P - Q)Q^{n-1}\rmN_{\epsilon} \right) \\ &\qquad + \sum_{j=1}^{n-2} \zeta_{1}( P^{n-j-1} - Q^{n-j-1} ) \, \nu_0\left( (P -Q )Q^{j}\rmN_{\epsilon} \right)  \\
	    &  
		\, \leq \,  \beta\epsilon + \zeta_{1}( P - \mathrm{N}) + (n-1) \frac{\alpha^{}_{\delta} \zeta^{}_{2,\delta}( P - \mathrm{N} )}{(n-1 + \epsilon^2)^{(1+\delta)/2}} \\ &\qquad + \sum_{j=1}^{n-2} \zeta_{1}( P^{n-j-1} - Q^{n-j-1} ) \frac{\gamma^{}_{\delta}
		\zeta^{}_{2,\delta}( P - \mathrm{N} )}{(j + \epsilon^2)^{(2+\delta)/2}} \\
		&  
		\, \leq \, \beta\epsilon + \zeta_{1}( P - \rmN ) + (n-1) \frac{\alpha^{}_{\delta} \zeta^{}_{2,\delta}(P - \rmN )}{(n-1 + \epsilon^2)^{(1+\delta)/2} } \\ &\qquad + \sum_{j=1}^{n-2} (n-j-1)^{(1-\delta)/2} \xi^{}_{\delta, 0} \frac{\gamma^{}_{\delta} \zeta^{}_{2,\delta}(P - \rmN )}{(j + \epsilon^2)^{(2+\delta)/2} } \\
		& 
		 \, = \,  \beta\epsilon + \zeta_{1}( P - \rmN )
		 + \sqrt{n} \alpha^{}_{\delta} \frac{n-1}{(n-1 + \epsilon^2)^{(1+\delta)/2} n^{(1-\delta)/2}} \frac{\zeta_{2,\delta}( \widetilde{P} - \rmN )}{n^{\delta/2}} \\
		 &\qquad + \sqrt{n}\xi^{}_{\delta, 0} \gamma^{}_{\delta} \frac{\zeta_{2,\delta}( P - \rmN )}{n^{\delta/2}} \sum_{j=1}^{n-2} \Big( \frac{n-j-1}{n} \Big)^{(1-\delta)/2} \frac{(j + \epsilon^2)^{(1-\delta)/2}}{(j + \epsilon^2)^{3/2}} \\
		&\, \leq \,  \beta\epsilon + \zeta_{1}( P - \rmN )+ \sqrt{n}  \alpha^{}_{\delta} \frac{\zeta_{2,\delta}(P - \rmN )}{n^{\delta/2}} \\
		&\qquad + \sqrt{n}\xi^{}_{\delta, 0} \gamma^{}_{\delta} \frac{\zeta_{2,\delta}( P - \rmN )}{n^{\delta/2}} \left(n +\epsilon^2\right)^{(1-\delta)/2} \sum_{j=1}^{\infty} \frac{1}{(j + \epsilon^2)^{3/2}}
	\end{align*}
	by using in the second step the smoothing inequality for $\zeta_{1}$, see e.g.~\cite[p. 86, Lemma~5.5]{Mattner2024},
	in the third the commutative ring identity~\eqref{Eq--ring_identity} multiplied
	by~$\rmN_{\epsilon}$ and the triangle inequality for $\zeta_{1}$,
	in the fourth the regularity~\eqref{Eq:regularity} of $\zeta_{1}$,
	  once with $M_1 = P-Q$, $M_2 = \rmN_{\epsilon}P^{n-1}$
	  and once with $M_1 = P^{n-j-1} - Q^{n-j-1}$, $M_2 = (P -Q )Q^{j}\rmN_{\epsilon}$,
	in the fifth $Q=\mathrm{N}$ and
	Lemma~\ref{Lem--(MNepsilon)_to_zeta2,delta(M)}
	with each $\sigma >1$ due to w.l.o.g.~$n\ge2$,
	and with  $C_{1, \delta} = (\frac{16\rme^{-1/2} - 4}{\sqrt{2\pi}})^{1-\delta}
	(\frac{4\rme^{-1/2}}{\sqrt{2\pi}})^{\delta} = \alpha^{}_{\delta}$, $C_{0, \delta} = (\frac{8\rme^{-1/2}}{\sqrt{2\pi}})^{1-\delta}
	(\frac{2 + 8\rme^{-3/2}}{\sqrt{2\pi}})^{\delta} = \gamma^{}_{\delta}$,
	in the sixth the inductive hypothesis~\eqref{Eq--Z/S_inductive_hypothesis}
	 to get, using the scaling identity~\eqref{Eq:scaling_identity_zeta_s} again,
	 \begin{align*}
		\zeta_{1}( P^{n-j-1} - Q^{n-j-1} )
		\, = \,  \sqrt{n-j-1} \zeta_{1}( \widetilde{P^{n-j-1}} - \rmN )
		\, \leq \,  \sqrt{n-j-1} \frac{\xi^{}_{\delta, 0}}{(n-j-1)_{}^{\delta/2}} \,,
	 \end{align*}
	 and in the final seventh and eighth steps nothing special.

	 We have
	 $\sum_{j=1}^\infty(j + \epsilon^2)^{-3/2}
	  < \sum_{j=1}^\infty\int_{j-1}^{j} (x+\epsilon^2)^{-3/2} \dd x
	  = \frac{2}{\epsilon}\,$.
	 Hence applying the above to $\epsilon = n^{(1-\delta)/2} \eta$ yields
	\begin{align*}
		\zeta_{1}( \widetilde{P^{\ast n}} - \rmN )
		&\, < \,  \frac{\zeta_{1}( P - \rmN )}{n^{1/2}}
		  \,+\, \alpha^{}_{\delta} \frac{\zeta_{2,\delta}( P - \rmN )}{n^{\delta/2}}
		  \,+\, \beta \frac{\eta}{n^{\delta/2}} \\
		  &\qquad \,+\, \xi^{}_{\delta, 0} \gamma^{}_{\delta}  \big( n( 1 + n^{-\delta}\eta^2 ) \big)^{(1-\delta)/2} \frac{2}{n^{(1-\delta)/2}\eta} \frac{\zeta_{2,\delta}( P - \rmN )}{n^{\delta/2}} \\
		&\, = \, \frac{\zeta_{1}( P - \rmN )}{n^{1/2}}
		  \,+\, \alpha^{}_{\delta} \frac{\zeta_{2,\delta}(P - \rmN )}{n^{\delta/2}}
		  \,+\, \beta \frac{\eta}{n^{\delta/2}} \\
		  &\qquad \,+\, \xi^{}_{\delta, 0} \gamma^{}_{\delta} 2 \frac{( 1 + n^{-\delta}\eta^2 )^{(1-\delta)/2}}{\eta} \frac{\zeta_{2,\delta}( P - \rmN )}{n^{\delta/2}} \\
		&\, \leq \,  \frac{1}{n^{\delta/2}_{}}  \Big( \zeta_{1}(P - \rmN )
		  \,+\, \alpha^{}_{\delta} \zeta^{}_{2,\delta}( P - \rmN )
		  \,+\, \beta \eta \Big.
		  \\&\qquad\qquad\qquad \Big. \,+\, \xi^{}_{\delta, 0} \gamma^{}_{\delta} 2 \frac{ (1 + \eta^2)^{(1-\delta)/2}}{\eta} \zeta_{2,\delta}( P - \rmN ) \Big) \\
		&\, \eqqcolon \,  \frac{1}{n^{\delta/2}_{}} \Big( A(\eta) \,+\, B(\eta) \xi_{\delta, 0} \Big)
	\end{align*}
	for $\eta \in \mathopen] 0, \infty \mathclose[\,$.
	Putting $\Eta \coloneqq \{\eta \in \mathopen] 0, \infty \mathclose[ : B(\eta) < 1\}$,
	we have	$\xi_{\delta, 0} = \inf_{\eta' \in \Eta} \frac{A(\eta')}{1-B(\eta')}$
	by the definition of $\xi^{}_{\delta}$ in~\eqref{Eq--Z/S_zeta1_zeta2,delta_1stPart_xi}.
	Hence
	\begin{align*}
		\zeta_{1}( \widetilde{P^{\ast n}} - \rmN )  &\, \leq \,  \frac{1}{n^{\delta/2}_{}} \inf_{\eta, \eta' \in \Eta} \Big( A(\eta) + B(\eta) \frac{A(\eta')}{1- B(\eta')} \Big)	\\
		&\, \leq \,  \frac{1}{n^{\delta/2}_{}} \inf_{\eta \in \Eta}\Big( A(\eta) + B(\eta) \frac{A(\eta)}{1 - B(\eta)} \Big)  \, = \,  \xi_{\delta, 0} 
	\end{align*}
	and hence~\eqref{Eq--Z/S_inductive_goal}.
    This completion of our inductive step
	proves~\ref{Thm--Z/S_zeta1_zeta2,delta_1stPart}.
	
	\smallskip
	\ref{Thm--Z/S_zeta1_zeta2,delta_2ndPart}
	With the notation from part~\ref{Thm--Z/S_zeta1_zeta2,delta_1stPart}
	 we have $\xi^{}_{\delta} \leq \xi^{}_{0}$ for $\delta \in \mathopen[0, 1 \mathclose]$,
	 due to $\alpha^{}_\delta\le \alpha^{}_0 = \alpha$
	 and  $\gamma^{}_\delta \le\gamma^{}_{0} =\gamma\,$.
	Let $(\vartheta, c) \in \mathopen] 0, 1 \mathclose[ \times \mathopen] 0, \infty \mathclose[$ be such that $\frac{\vartheta^2}{4\gamma^2} - \frac{1}{c^2} > 0$ holds.
	For $\eta\coloneqq c\zeta$ with $\zeta \in \mathopen]0,\infty\mathclose[$
	and with recalling~\eqref{Eq--Z/S_zeta1_zeta2,delta_1stPart_g},
	we have the equivalence
	\begin{align*}
		\gamma g^{}_{0}(\eta) \zeta
		 \, = \,  2 \gamma \frac{\sqrt{ 1+ \eta^2}}{\eta }  \zeta
		   \, = \, 2 \gamma \sqrt{\zeta^2 + \frac{1}{c^2}}
		    \, \leq \, \vartheta
		    \,\, \Leftrightarrow \,\,  \zeta^2	+ \frac{1}{c^2}
		    \, \leq \,  \frac{\vartheta^2}{4 \gamma^2} \,,
	\end{align*}
	yielding for $\delta \in \mathopen[ 0 , 1 \mathclose]$ and $P \in \cP_{2+\delta}$ satisfying
	$\zeta_{2,\delta}(\widetilde{P} - \rmN) \leq \sqrt{\frac{\vartheta^2}{4\gamma^2} - \frac{1}{c^2}}$
	the boundings
	\begin{align*}
		\xi^{}_{\delta} \left( \zeta_{1}(\widetilde{P} - \rmN), \zeta_{2,\delta}(\widetilde{P} - \rmN) \right)  \, &\leq \,\xi^{}_{0} \left( \zeta_{1}(\widetilde{P} - \rmN), \zeta_{2,\delta}(\widetilde{P} - \rmN) \right)  \\ &\leq \, \frac{\zeta_{1} + \alpha\zeta_{2,\delta} + \beta c \zeta_{2,\delta}}{1 - \vartheta}(\widetilde{P} - \rmN)  \\ &\leq \,     \frac{1 + \alpha + \beta c}{1 - \vartheta} \big( \zeta_{1} \vee \zeta_{2,\delta} \big)(\widetilde{P} - \rmN)
	\end{align*}
	and thus, using~\ref{Thm--Z/S_zeta1_zeta2,delta_1stPart},
	inequality~\eqref{Eq--Z/S_zeta1_zeta2,delta_2ndPart}.
\end{proof}
\begin{proof}[Proof of Theorem~\ref{Thm:Lyapunov-Katz_for_summands_Z-close_to_normal}]
	Let $\delta\in[0,1]$,
	$P\in\cP_2$, and $2\le n\in\N$.
	
	\smallskip 1. Let $\Xi:\cP_2\rightarrow[0,\infty]$ be a functional such that we have
	\begin{align} \label{Eq--PoT1.1_functional}
		\zeta_1\big(\widetilde{P^{\ast k}}-\mathrm{N}\big)
		&\,\le\, \frac{\Xi(P)}{k^{\delta/2}} \quad\text{ for $P\in\cP_2$ and $k\in\N$\,}.
	\end{align}
	Since $n=2k+\rho$ with $k\in\N$ and $\rho\in\{0,1\}$,
	an application of the convolution inequality as in \cite[p.~72, proof of Theorem~1.5]{Mattner2024}
	then yields
	\begin{align*}
		\big\| \widetilde{P^{\ast n}}-\mathrm{N} \big\|^{}_\mathrm{K}
		&\,\le\, \frac{\Xi(P)}{2\sqrt{2\pi} n^{\delta/2}}
		\left( \left(\frac{(2k+\rho)^\delta k^{1-\delta}}{k+\rho}\right)^{\!\frac{1}{4}}
		+ \left(\frac{(2k+\rho)^\delta (k+\rho)^{1-\delta}}{k}\right)^{\!\frac{1}{4}}  \right)^2,
	\end{align*}
	and here both numerators of the fractions in parentheses can be upper bounded
	by replacing $\delta$ with $1$.
	Hence, by \cite[pp.~72-73, the calculations leading to (3.16)]{Mattner2024}, we get
	\begin{align} \label{Eq--PoT1.1_conv_yields}
		\big\| \widetilde{P^{\ast n}}-\mathrm{N} \big\|^{}_\mathrm{K}
		&\,\le\, c_1\frac{\Xi(P)}{n^{\delta/2}}
		\quad\text{ with } c_1\,\coloneqq\,
		\frac{\sqrt{3}(2^{-1/4}+1)^2}{2\sqrt{2\pi}}
		\, =\, 1.1708\ldots\,.
	\end{align}
	
	\smallskip2. Theorem~\ref{Thm--Z/S_zeta1_zeta2,delta}\ref{Thm--Z/S_zeta1_zeta2,delta_1stPart} with $\xi^{}_{\delta} \leq \xi^{}_{0}$ yields~\eqref{Eq--PoT1.1_functional} with $\Xi(P) \coloneqq \xi^{}_{0}( \zeta_{1}(\widetilde{P} - \rmN), \zeta_{2,\delta}(\widetilde{P} - \rmN) )$
	where $\xi^{}_{0}$ satisfies
	\begin{align} \label{Eq--PoT1.1_xi0}
	 \xi^{}_{0}(\varkappa, \zeta)
		\, \leq \, \frac{\varkappa + \alpha\zeta + \beta \eta}
			{ \left(1 - 2\gamma\sqrt{ 1 + \eta^{-2}}  \, \zeta \right)_{+} }
		\quad \text{for $(\varkappa, \zeta) \in \mathopen[ 0, \infty \mathclose[^2$ and $\eta \in \mathopen] 0, \infty \mathclose[$} \,,
	\end{align}
	hence~\eqref{Eq--PoT1.1_conv_yields} yields
	\begin{align} \label{Eq--PoT1.1_Estimate_for_zeta_small}
	 \big\| \widetilde{P^{\ast n}}-\mathrm{N} \big\|^{}_\mathrm{K}
		\, \leq \,
		 \frac{c_1}{n^{\delta/2}}
		 \frac{\zeta_{1} + \alpha\zeta_{2,\delta} + \beta \lambda\zeta_{2,\delta} }
		  { \left(1 - 2\gamma\sqrt{ \zeta^2_{2,\delta}
		   + \lambda_{\phantom{2}}^{-2} } \, \right)_{+} }
		  ( \widetilde{P} - \rmN ) \quad \text{for $\lambda \in \mathopen] 0, \infty \mathclose[$}
	\end{align}
	by choosing $\eta = \lambda\zeta$ in~\eqref{Eq--PoT1.1_xi0}.
	
	Inequality~\eqref{Eq:Lyapunov-Katz} with $c = 1.8546 \eqqcolon c_{2}$ from
	Korolev and Dorofeeva~(2017)
\cite[p.~46, Corollary~3, definition of $C_1$ on page 39]{KorolevDorofeeva2017}
	and Lemma~\ref{Lem--nu2+delta_to_zeta2,delta} yield
	\begin{align} \label{Eq--PoT1.1_Estimate_for_zeta_large}
		\big\| \widetilde{P^{\ast n}}-\mathrm{N} \big\|^{}_\mathrm{K}  \, \leq \,  c_{2} \frac{\nu_{2+\delta}(\widetilde{P})}{n^{\delta/2}}  \, \leq \,  \frac{c_{2}}{n^{\delta/2}} \left( 6 \zeta_{2,\delta}(\widetilde{P} - \rmN) + \omega \right) \,,
	\end{align}
	with $\omega = \frac{25}{24}\frac{4}{\sqrt{2\pi}} + \frac{4}{27}$, and in combination with~\eqref{Eq--PoT1.1_Estimate_for_zeta_small} we then have
	\begin{align*}
		\big\| \widetilde{P^{\ast n}} - \rmN \big\|^{}_{\rmK}  \, \leq \,  \frac{c(\lambda)}{n^{\delta/2}} \left( \zeta_{1} \vee \zeta_{2,\delta} \right)(\widetilde{P} - \rmN) \quad \text{for $\lambda \in \mathopen] 0, \infty \mathclose[$} \,,
	\end{align*}
	with
	\begin{align*}
		c(\lambda)  \, &\coloneqq \,  \sup_{\varkappa, \zeta > 0}
		  \frac{1}{\varkappa \vee \zeta} \left(
		   \left(
		     c_1 \frac{\varkappa + \alpha\zeta + \beta \lambda \zeta}
		         {\: 1 - 2\gamma\sqrt{ \zeta^2 + \lambda^{-2}  }  \:}
		    \right)_{\!\!+}
		   \,\wedge\: \Big( c_{2} \left(6\zeta + \omega \right)\Big) \right) \\
		&= \, \sup_{\zeta>0} \left(
		 \left(
		 c_1 \frac{1 + \alpha + \beta \lambda }{\:1 - 2\gamma\sqrt{ \zeta^2
		 + \lambda^{-2}} \:}
		 \right)_{\!\!+}
		\,\wedge\: \left(c_{2} \left(6 + \frac{\omega}{\zeta} \right) \right)\right)
		\quad \text{for $\lambda \in \mathopen] 0, \infty\mathclose[$} \,.
	\end{align*}
	If $\lambda > 2\gamma = 3.87153\ldots$,
	then $2\gamma\sqrt{\lambda^{-2}} < 1$,	and by monotonicity in $\zeta$
	then $c(\lambda)$ is finite and attained uniquely
	at the positive solution $\zeta^{\ast}(\lambda)$ to the fourth-degree equation
	\begin{align*}
		&\frac{A}{1-B\sqrt{\zeta^2 + C}} \, = \,  D + \frac{E}{\zeta}
	\end{align*}
	with $A \coloneqq c_1(1+\alpha+\beta\lambda)/c_2\,$, $B \coloneqq 2 \gamma\,$,
	$C \coloneqq \lambda^{-2}$, $D \coloneqq 6\,$, $E \coloneqq \omega\,$.

	Hence the claim~\eqref{Eq:Lyapunov-Katz-Senatov_for_summands_Z-close_to_normal} holds with
	$c \coloneqq 	\inf_{\lambda >2\gamma}	c(\lambda) < \infty$,
	which should be computed numerically and can be upper bounded by $c(8.5) = 47.10171\ldots<48$.
\end{proof} 

\section{Auxiliary Results} 
We put, for use in the proofs of Lemmas \ref{Lem--Scaling_ineq_for_zeta2,delta}, \ref{Lem--(MNepsilon)_to_zeta2,delta(M)} and \ref{Lem--nu2+delta_to_zeta2,delta} below,
\begin{align*}
	&\cF_{2,\delta} \,\coloneqq\, \cF_{2, (\cdot) \wedge (\cdot)_{}^{\delta}} \,, \\
	&\cF_{2,\delta}^{\infty} \,\coloneqq\, \cF_{2, (\cdot) \wedge (\cdot)^{\delta}}^{\infty} \,=\, \left\{ f \in \cF_{2,\delta} \,:\; f \ \text{bounded} \right\}
\end{align*}
for $\delta \in \mathopen[0, 1\mathclose]$. We recall that $\cF_{2,\delta}^{\infty}$ appears in the definition of $\zeta_{2,\delta}$
in \eqref{Eq:Def_zeta_2,delta} and \eqref{Eq:Def_zeta_m,g}.

Part~\ref{part:zeta_m,g_for_unbounded_f} of the following lemma states that the boundedness condition
on the functions~$f$ in the definition of $\zeta_{m,g}(M)$ may be omitted
for $M\in\cM_{m,m}\cap\cM_{m,g}$ if $g\in\cG$ is unbounded.
We suspect that this last condition is actually unnecessary, but our present proof uses it.
In part~\ref{part:primitive_of_g},
``primitive of order $m$'' refers to $m$ differentiations  on $]0,\infty[\,$,
but only
$(m-1)_+$ differentiations at zero.
\begin{lemma}   \label{Lem:cF_m,g-bounds_and_zeta_m,g_for_unbounded_f}
Let $m\in\N_0$ and $g\in\cG$.

\begin{parts}
\item  For the primitive $g^{(-m)}$ of $g$ of order $m$ and \label{part:primitive_of_g}
 with derivatives of all orders $\le$~$m-1$ vanishing at zero, we have
\begin{equation}                 \label{Eq:mth_primitve_of_g_represented_and_bounded}
  g^{(-m)}(u) \,=\,
   \left\{\!
    \begin{array}{ll}
     g (u)                                            &\text{if }m=0, \\
     \int_0^{u}g(y)\frac{(u-y)^{m-1}}{(m-1)!}\dd y &\text{if }m\ge1
    \end{array}
    \!\right\}
    \le \, g(u)\frac{u^m}{m!}
    \quad\text{ for } u\in [0,\infty[\,.
 \end{equation}
\item Let $f\in\cF_{m,g}$ and $x\in\R$. \label{part:Bounds_for_cF_m,g-functions}
 Then the Taylor remainder
 $f_0(x) \coloneqq f(x)-\sum_{j=0}^m \frac{f^{(j)}(0)}{j!}x^j$
 of order $m$ satisfies
 $|f_0(x)| \,\le\, g^{(-m)}\big(|x|\big)\,\le\, g\big(|x|\big)\frac{|x|^m}{m!}$.
\item $g^{(-m)}\circ|\cdot|\in\cF_{m,g}$
 iff  $m$ is even or $g$ is linear.                   \label{part:When_primitive of_g_in_cF_m,g}
\item                                                  \label{part:zeta_le_nu}
 $\zeta_{m,g}(M) \,\le\, \nu_{0,g^{(-m)}}(M)  \,\le\, \frac{1}{m!}\nu_{m,g}(M)$ for $M\in\cM_{m,m}$.
\item                                             \label{part:zeta_s_le_nu_s}
 $\zeta_s(M) \,\le\,\left(\prod_{j=1}^m (j+\delta)\right)^{-1}\nu_s(M)$
 for $s\in\mathopen]0,\infty\mathclose[\,$, $M\in\cM_{s,m}\,$, $\delta\coloneqq s-m\in[0,1]$.
\item If $g$ is unbounded, then      \label{part:zeta_m,g_for_unbounded_f}
 $\zeta_{m,g}(M)=\sup_{f\in\cF_{m,g}}\left|\int f\,\dd M\right|$
 for $M\in\cM_{m,m}\cap \cM_{m,g}$.
\end{parts}
\end{lemma}
\begin{proof} We may assume that $g(0)=g(0+)$, so that then $g$ is continuous everywhere.

\smallskip\ref{part:primitive_of_g} Obvious by continuity and isotonicity of $g$.

\smallskip\ref{part:Bounds_for_cF_m,g-functions}
For $f\in\cF_{m,g}$ the condition \eqref{Eq:Increment_bound_for_mth_derivative_with_g}
with $y=0$ reads
\begin{equation} \label{Eq:Increment_bound_for_mth_derivative_with_g_and_y=0}
 \left|f^{(m)}(x)- f^{(m)}(0)\right| \,\le \, g\big(|x|\big)\quad\text{ for }x\in\R\,.
\end{equation}
If $m=0$, then \eqref{Eq:Increment_bound_for_mth_derivative_with_g_and_y=0}
directly yields the claim.  
If $m\ge 1$, then the integral formula for the Taylor remainder of
order~$m-1$ yields
\begin{align*}
 f_0(x)&\,=\, -\frac{f^{(0)}(0)}{m!}x^m
    + \int_0^x f^{(m)}(y)\frac{(x-y)^{m-1}}{(m-1)!}\dd y \\
  &\,=\,\int_0^x \left(f^{(m)}(y)-f^{(m)}(0)\right)\frac{(x-y)^{m-1}}{(m-1)!}\dd y
\end{align*}
and hence the claim via~\eqref{Eq:Increment_bound_for_mth_derivative_with_g_and_y=0}
and~\eqref{Eq:mth_primitve_of_g_represented_and_bounded}.

\smallskip\ref{part:When_primitive of_g_in_cF_m,g}
The subadditivity of $g$, noted on page~\pageref{page:g_subadditive}
after~\eqref{Eq:Def_g_b}, and the isotonicity of $g$ yield
$g\big(|x|\big)\le g\big(|x-y| +|y|\big) \le g\big(|x-y|\big) +g\big(|y|\big)$
for $x,y\in\R$, and hence
\begin{equation}    \label{Eq:increment_of_g_bounded_by_g_of_increment}
 \big| g\big(|x|\big) -g\big(|y|\big)\big|
 \,\le\, g\big(|x-y|\big) \quad\text{ for }x,y\in\R\,.
\end{equation}

Let us put $f_{m,g}\coloneqq g^{(-m)}\circ|\cdot|$.
If $m=0$, then \eqref{Eq:mth_primitve_of_g_represented_and_bounded} and
\eqref{Eq:increment_of_g_bounded_by_g_of_increment}
yield $f_{m,g}= g\circ|\cdot|\in\cF_{m,g}\,$.

Let now $m\ge1$.
The function $f_{m,g}$ is $m-1$ times differentiable
with $f^{(m-1)}_{m,g}(x)=\sgn(x)^{m-1}\int_0^{|x|}g(y)\,\dd y$ for $x\in\R$,
and with $m$th derivative $f^{(m)}_{m,g}(x)=\sgn(x)^m g\big(|x|\big)$
for $x\in\R\setminus\{0\}$.
Hence $f_{m,g}$ is $m$ times differentiable iff $m$ is even
or $g(0)=0$.
If $m$ is even, then
we get~\eqref{Eq:Increment_bound_for_mth_derivative_with_g} for $f=f_{m,g}$
since
$\big| \sgn(x)^mg\big(|x|\big) -\sgn(y)^mg\big(|y|\big)\big|
 = \big| g\big(|x|\big)-g\big(|y|\big)\big| $ for $x,y\in\R$;
hence then  $f_{m,g}\in\cF_{m,g}\,$.

If however $m$ is odd and $f_{m,g}\in\cF_{m,g}\,$,
then $g(0)=0$, and
for $u,v\in\mathopen[0,\infty[$ and with $x\coloneqq u$ and $y\coloneqq -v$
we get
$g(u)+g(v) = \big| \sgn(x)^mg\big(|x|\big) -\sgn(y)^mg\big(|y|\big)\big|
= \big|f_{m,g}^{(m)}(x)- f_{m,g}^{(m)}(y)\big|
\le g\big(|x-y|\big) =g(u+v)$,
so that $g$ is superadditive, hence additive by subadditivity,
hence linear by continuity.
That linearity of $g$ implies $f_{m,g}\in\cF_{m,g}\,$ is easy to see.

\smallskip\ref{part:zeta_le_nu} If $M\in\cM_{m,m}$ and $f\in\cF^\infty_{m,g}$,
then with $f_0$ as in part~\ref{part:Bounds_for_cF_m,g-functions}
we have $\int f\,\dd M = \int f_0\,\dd M$,
and by~\ref{part:Bounds_for_cF_m,g-functions} hence
$\left| \int f\,\dd M \right| \,\le\, \int \left|f_0\right|\,\dd|M|
\,\le\, \nu_{0,g^{(-m)}}(M) \,\le\, \frac{1}{m!}\nu_{m,g}(M)$.

\smallskip\ref{part:zeta_s_le_nu_s} Follows from part~\ref{part:zeta_le_nu}:
If $\delta>0$, let there $g\coloneqq(\cdot)^\delta$.
If $\delta=0$, then $m\ge1$, and we can apply what has just been proved
to $(m-1,1)$ in place of $(m,\delta)$.

\smallskip\ref{part:zeta_m,g_for_unbounded_f} Let $M\in\cM_{m,m}\cap \cM_{m,g}$.
By isotonicity and unboundedness of $g$, we have
$\lim_{u\rightarrow\infty}g(u)=\infty$.

Let $f\in\cF_{m,g}$.
We then have, with $f_0$ as in part~\ref{part:Bounds_for_cF_m,g-functions},
a sequence $(f_n)$ in $\cF_{m,g}^\infty$ with
$f_n\rightarrow f_0$ pointwise and
$|f_n(x)|\le g\big(|x|)\frac{|x|^m}{m!}$ for $n\in\N$ and $x\in\R$.
This claim generalises, for the case of functions with a vanishing Taylor polynomial
of order $m$ around zero,
the main part of \cite[p.~504, Lemma 2.2]{MattnerShevtsova2019},
which considers the special case of $g=(\cdot)^\alpha$ with $\alpha\in\mathopen]0,1\mathclose]$.
The present claim can be proved as its special case, replacing $w^\alpha$ at every occurrence
with $g(w)$ with the present $g$, and using~\ref{part:Bounds_for_cF_m,g-functions}
in place of the argument starting with \cite[p.~504, (2.8)]{MattnerShevtsova2019}.
\pr{The unboundedness of $g$ is used here for ensuring that
 the present analogue of $c\|y\|(b_n-n)^{-\alpha}$ in
 \cite[p.~504, fifth line of the proof of Lemma~2.2]{MattnerShevtsova2019},
 namely $c\|y\|/g(b_n-n)$, can be made $\le$~$\frac{1}{2n}$ by choosing
 $b_n$ large enough.}
Hence $\int f\,\dd M = \int f_0\,\dd M = \lim_{n\rightarrow\infty} \int f_n\,\dd M$
by dominated convergence, and hence $ \left|\int f\,\dd M\right| \le  \zeta_{m,g}(M)$.

This proves $\sup_{f\in\cF_{m,g}}\left|\int f\,\dd M\right|\le \zeta_{m,g}(M)$,
and hence equality.
\end{proof}

\begin{lemma}[A scaling inequality for $\zeta_{2,\delta}$] \label{Lem--Scaling_ineq_for_zeta2,delta}
	Let $\delta \in \mathopen[ 0, 1 \mathclose]$ and $M \in \cM$. We have
	\begin{align} \label{Eq:Scaling_ineq_for_zeta_2,delta}
		\zeta_{2,\delta}\big( M(\tfrac{\cdot}{a}) \big)  
		\, \leq \,  \big(a^{3} \, \vee \, a^{2+\delta} \big) \zeta_{2,\delta}( M ) \quad \text{for $a \in \mathopen] 0, \infty \mathclose[$} \,.
	\end{align}
\end{lemma}
\begin{proof}
	Let $a \in \mathopen] 0, \infty \mathclose[$. For $f \in \cF_{2,\delta}^{\infty}$ we put $f_a \coloneqq f(a\cdot)$. Then $f_a$ is bounded and twice differentiable with $f_a'' = a^2 f''(a\cdot)$, and for any $x,y \in \R$ we have
	\begin{align*}
		\abs{ f_a''(x) - f_a''(y) }  &\, = \,  a^2 \abs{ f''(ax) - f''(ay) }  \\
		&\, \leq \,  a^2 \big( \abs{ ax - ay }^{\delta} \wedge \abs{ ax - ay } \big)  \, \leq \,  (a^3 \vee a^{2+\delta}) \big( \abs{ x-y }^{\delta} \wedge \abs{ x-y } \big) \,.
	\end{align*}
	Therefore, $(a^{3} \vee a^{2+\delta})^{-1} f_a \in \cF_{2,\delta}^{\infty}$ and we can conclude that
	\begin{align*}
		\zeta_{2,\delta}\big( M(\tfrac{\cdot}{a}) \big)
		&\, = \,  \sup_{f \in \cF^{\infty}_{2,\delta}} \Big| \int f_a(x) \,\dd M(x) \Big| \\
		&\, = \,  \big( a^3 \vee a^{2+\delta} \big) \sup_{f \in \cF^{\infty}_{2,\delta}} \Big| \int \big( a^3 \vee a^{2+\delta} \big)^{-1} f_a(x) \,\dd M(x) \Big| \\
		&\, \leq \,  \big( a^3 \vee a^{2+\delta} \big) \sup_{f \in \cF^{\infty}_{2,\delta}} \Big| \int f(x) \,\dd M(x) \Big| \\
		&\, = \,  \big( a^3 \vee a^{2+\delta} \big) \zeta_{2,\delta}( M ) \,. 
	\end{align*}
\end{proof}

\begin{lemma} \label{Lem--(MNepsilon)_to_zeta2,delta(M)}
	With $\phi$ denoting the standard normal density let
	\begin{align*}
		D_{2}  \, &\coloneqq \,  \int \abs{ \phi''(z) } \,\dd z  \, = \,  \frac{4\rme^{-1/2}}{\sqrt{2\pi}}  \, = \,  0.967882\ldots \,, \\
		D_{2,1}  \, &\coloneqq \,  \int \abs{z} \abs{ \phi''(z) } \,\dd z  \, = \,  \frac{8\rme^{-1/2} - 2}{\sqrt{2\pi}}  \, = \,  1.13788\ldots \,, \\
		D_{3}  \, &\coloneqq \,  \int \abs{ \phi'''(z) } \,\dd z  \, = \,  \frac{2 + 8\rme^{-3/2}}{\sqrt{2\pi}}  \, = \,  1.510013\ldots \,.
	\end{align*}
	Then for $\delta \in \mathopen[ 0, 1 \mathclose]$ the following holds: 
	\begin{parts}
	\item \label{Lem--nuo(MNepsilon)_to_zeta2,delta(M)}
	There exists a constant $C_{0, \delta} \in \mathopen] 0, \infty \mathclose[$ bounded in $\delta$ and satisfying
	\begin{align*}
		\nu_{0}( M\rmN_{\sigma} )  \, \leq \,  C_{0, \delta} \frac{\zeta_{2,\delta}(M )}{\sigma^{3} \wedge \sigma^{2+\delta}} \quad \text{for $\sigma \in \mathopen] 0, \infty[$ and $M \in \cM$} \,. 
	\end{align*}
	We may take $C_{0, \delta} = (2D_{2})^{1-\delta}  D_{3}^{\delta}$ with $C_{0, \delta} \leq C_{0} \coloneqq 2D_{2} = 1.935765\ldots$\;.
	
	\item \label{Lem--zeta1(MNepsilon)_to_zeta2,delta(M)} 	
	There exists a constant $C_{1, \delta} \in \mathopen] 0, \infty \mathclose[$ bounded in $\delta$ and satisfying
	\begin{align*} 
		\zeta_{1}( M\rmN_{\sigma} )  \, \leq \,  C_{1, \delta} \frac{\zeta_{2,\delta}( M )}{\sigma^{2} \wedge \sigma^{1+\delta}} \quad \text{for $\sigma \in \mathopen] 0, \infty[$ and $M \in \cM$} \,. 
	\end{align*}
	We may take $C_{1, \delta} = (2D_{2,1})^{1-\delta}  D_{2}^{\delta}$ with $C_{1, \delta} \leq C_{1} \coloneqq 2D_{2,1} = 2.27576\ldots$\;.
	\end{parts}
\end{lemma}
\begin{proof}
	The stated values of $D_{2}, D_{2,1}, D_{3}$ are easily checked. 
	
	Part~\ref{Lem--nuo(MNepsilon)_to_zeta2,delta(M)} is, in the case of $\sigma = 1$ and $M = P-Q$ for probability measures $P$ and $Q$, and without an explicit value for $C_{0}$, contained in Senatov~(1980) \cite[Lemma~1]{Senatov1980}. In deriving an explicit value for $C_{0}$ and in proving part~\ref{Lem--zeta1(MNepsilon)_to_zeta2,delta(M)}, we follow Mattner and Shevtsova~(2019) \cite[pp. 513-515]{MattnerShevtsova2019} who in turn follow the outline of the reasoning employed in Senatov~(1998) \cite[Lemma~2.10.1]{Senatov1998}. As in \cite[p.~57, (1.61, 1.64)]{Mattner2024}
	and in accordance with \cite[p.~497]{MattnerShevtsova2019}, we let $\cF^\infty_0$ denote the set of all
	(bounded) measurable functions $f:\R\rightarrow\R$ with $\sup_{x\in\R}|f(x)|\le 1$,
	and $\cF_r^\infty\coloneqq \cF^\infty_{r-1,(\cdot)}$ for $r\in\N$, so that $\cF_1^\infty$ is the set of all bounded functions
	with Lipschitz constant $\le$~$1$.
	 
	Let $\sigma \in \mathopen] 0, \infty\mathclose[$\,, $\delta \in \mathopen[ 0, 1 \mathclose]$, and $M \in \cM$.
	
	For $f:\R\rightarrow\R$ bounded and measurable, we
	define a bounded function $h$
	and  differentiate under the integral:
	\begin{align}
	  h(x)      &\,\coloneqq\,  \int f(x+z)\phi(z) \,\dd z
	              \, = \,  \int f(z)\phi(x- z) \,\dd z\,, \label{Eq:h_in_proof_L_3.3}  \\
	  h^{(k)}(x)& \, = \,  \int f(z) \phi^{(k)}(x-z) \,\dd z
	              \, =\,  \int f(x+z) \phi^{(k)}(z) \,\dd z
	                                                  \label{Eq:h^(k)_in_proof_L_3.3}
	\end{align}
    for $x\in\R$ and $k\in\N_0$.

	\smallskip
	\ref{Lem--nuo(MNepsilon)_to_zeta2,delta(M)} We may assume w.l.o.g.~$\sigma = 1$, as for any $\sigma \in \mathopen] 0, \infty \mathclose[$ we have, given the validity of~\ref{Lem--nuo(MNepsilon)_to_zeta2,delta(M)} for $\sigma = 1$ and using the scaling inequality~(\ref{Eq:Scaling_ineq_for_zeta_2,delta}) for $\zeta_{2,\delta}$,
	\begin{align*}
		\nu_{0}( M\rmN_{\sigma} )  
		&\, = \,  \nu_{0}\big( (M(\sigma\cdot)\rmN) (\tfrac{\cdot}{\sigma}) \big) \\
		&\, = \, \nu_{0}\big( M(\sigma\cdot)\rmN \big) \\ 
		&\, \leq \,  C_{0, \delta} \, \zeta_{2,\delta}\big( M(\sigma\cdot) \big)  \\
		&\, \leq \,  C_{0, \delta} \big( (\tfrac{1}{\sigma})^{3} \vee (\tfrac{1}{\sigma})^{2+\delta} \big) \zeta_{2,\delta}( M )  \\
		&\, = \,  C_{0, \delta} \frac{\zeta_{2,\delta}( M )}{\sigma^{3} \wedge \sigma^{2+\delta}} \,.
	\end{align*}
	For $f \in \cF_{0}^{\infty}$
	and with $h$ from~\eqref{Eq:h_in_proof_L_3.3} we then get
	from~\eqref{Eq:h^(k)_in_proof_L_3.3}
	for the absolute difference $\abs{h''(x) - h''(y)}$ with $x,y \in \R$ the two bounds
	\begin{align} \label{Eq--nu0_h_absolute}
		\abs{ h''(x) - h''(y) }  \, \leq \, 2 \norm{h''}_{\infty}  \, \leq \,  2 \int \abs{ \phi''(z) } \,\dd z  \, = \,  2 D_{2}
	\end{align}
	and
	\begin{align} \label{Eq--nu0_h_difference}
		\abs{h''(x) - h''(y) }  \, \leq \, \norm{h'''}_{\infty} \cdot \abs{ x - y}  \, \leq \,  \int \abs{ \phi'''(z) } \,\dd z \cdot \abs{ x - y }  \, = \,  D_{3} \abs{ x - y } \,. 
	\end{align}
	Taking the geometric mean of~\eqref{Eq--nu0_h_absolute} and~\eqref{Eq--nu0_h_difference} with the exponents $1-\delta$ and $\delta$ yields
	\begin{align} \label{Eq--nu0_h_geometric_mean}
		\abs{ h''(x) - h''(y) }  \, \leq \,  (2D_{2})^{1-\delta} \cdot D_{3}^{\delta} \cdot \abs{ x - y }^{\delta}
	\end{align}
	which in turn yields, combined with~\eqref{Eq--nu0_h_difference},
	\begin{align*}
		\abs{ h''(x) - h''(y) }  \, \leq \, C_{0 , \delta}\,
		\big( \abs{ x - y }^{\delta} \wedge \abs{ x - y } \big)
	\end{align*}
	where
	 $C_{0 , \delta} \coloneqq \big( (2D_{2})^{1-\delta}  D_{3}^{\delta} \big)
	  \vee  D_{3}  = (2D_{2})^{1-\delta}  D_{3}^{\delta} \leq 2D_{2} \eqqcolon C_{0}$.
	 Therefore we have $\frac{1}{C_{0, \delta}} h \in \cF_{2,\delta}^{\infty}$ and can conclude~\ref{Lem--nuo(MNepsilon)_to_zeta2,delta(M)} with
	\begin{align} \label{Eq--nu0_conclusion(a)}
		\nu_{0}(M\rmN)  \, = \,  C_{0, \delta }\sup_{f \in \cF^{\infty}_{0}}  \left| \int \frac{1}{C_{0,\delta}} \int f(x+z)\phi(z) \,\dd z \,\dd M(x) \right|  \, \leq \,  C_{0, \delta} \, \zeta_{2,\delta}(M) \,.
	\end{align}
	
	\smallskip
	\ref{Lem--zeta1(MNepsilon)_to_zeta2,delta(M)} As in the first part, we may assume w.l.o.g.~$\sigma = 1$, as for any $\sigma \in \mathopen] 0, \infty [$ we have, given the validity of~\ref{Lem--zeta1(MNepsilon)_to_zeta2,delta(M)} for $\sigma = 1$ and using the homogeneity~\eqref{Eq:scaling_identity_zeta_s} of $\zeta_{1}$ and the scaling inequality~(\ref{Eq:Scaling_ineq_for_zeta_2,delta}) for $\zeta_{2,\delta}$,
	\begin{align*}
		\zeta_{1}( M\rmN_{\sigma} )  
		&\, = \,  \zeta_{1}\big( (M(\sigma\cdot)\rmN) (\tfrac{\cdot}{\sigma}) \big) \\
		&\, = \, \sigma \zeta_{1}\big( M(\sigma\cdot)\rmN \big) \\ 
		&\, \leq \,  C_{1, \delta} \sigma \, \zeta_{2,\delta}\big( M(\sigma\cdot) \big)  \\
		&\, \leq \,  C_{1, \delta} \sigma \big( (\tfrac{1}{\sigma})^{3} \vee (\tfrac{1}{\sigma})^{2+\delta} \big) \zeta_{2,\delta}( M )  \\
		&\, = \,  C_{1, \delta} \frac{\zeta_{2,\delta}( M )}{\sigma^{2} \wedge \sigma^{1+\delta}} \,.
	\end{align*}
	For $f \in \cF_{1}^{\infty}$ 	and with $h$ from~\eqref{Eq:h_in_proof_L_3.3} we then get
	from~\eqref{Eq:h^(k)_in_proof_L_3.3},
    using in the second step below $\abs{ f(a) - f(b) } \leq \abs{ a - b}$ for $a, b \in \R$,
	\begin{align} \label{Eq--zeta1_h_difference}
		\abs{ h''(x) - h''(y) }  \, \leq \,  \int \abs{ f(x+z) - f(y+z) } \abs{ \phi''(z) } \,\dd z  \, \leq \, D_{2} \abs{ x - y } \,,
	\end{align}
	and, using in the first step below $\int \phi''(z) \,\dd z = 0$,
	\begin{align} \label{Eq--zeta1_h''_infty}
		\abs{ h''(x) }  \, = \,  \left| \int \big( f(x+z) - f(x) \big) \phi''(z) \,\dd z \right|  \, \leq \,  \int \abs{ z }  \left| \phi''(z) \right| \,\dd z  \, = \,  D_{2,1} \,, 
	\end{align}
	which then yields
	\begin{align} \label{Eq--zeta1_h_absolute}
		\abs{ h''(x) - h''(y) }  \, \leq \, 2 \norm{h''}_{\infty}  \, \leq \,  2 D_{2,1}
	\end{align}
	for $x,y \in \R$. Taking the geometric mean of~(\ref{Eq--zeta1_h_absolute}) and~(\ref{Eq--zeta1_h_difference}) with the exponents $1-\delta$ and $\delta$ we obtain
	\begin{align} \label{Eq--zeta1_h_geometric_mean}
		\abs{ h''(x) - h''(y) }  \, \leq \,  (2D_{2,1})^{1-\delta} \cdot D_{2}^{\delta} \cdot \abs{ x - y }^{\delta} \,,
	\end{align}
	which then, combined with~(\ref{Eq--zeta1_h_difference}) yields
	\begin{align*}
		\abs{ h''(x) - h''(y) }  \, \leq \, C_{1,\delta} \, \big( \abs{ x - y }^{\delta} \wedge \abs{ x - y } \big)
	\end{align*}
	where $C_{1, \delta} \coloneqq \big((2D_{2,1})^{1-\delta}  D_{2}^{\delta}\big) \vee D_{2} = (2D_{2,1})^{1-\delta}  D_{2}^{\delta} \leq 2D_{2,1} \eqqcolon C_{1}$. Therefore we have $\frac{1}{C_{1, \delta}} f \in \cF_{2,\delta}^{\infty}$ and analogously to \eqref{Eq--nu0_conclusion(a)} we can conclude~\ref{Lem--zeta1(MNepsilon)_to_zeta2,delta(M)}.
\end{proof}

\begin{lemma}                               \label{Lem--nu2+delta_to_zeta2,delta}
	Let $\delta \in \mathopen[ 0, 1 \mathclose]$. We have
	\begin{align*}
		\nu_{2+\delta}(P)  \, \leq \,  (1+\delta)(2+\delta) \, \zeta_{2,\delta}(P - \rmN) + \int h^{}_{\delta} \,\dd \rmN \quad
		\text{for $P \in \widetilde{\cP_{2+\delta}}$}\,,
	\end{align*}
	where $h^{}_{\delta}$ is the function defined by~\eqref{Eq--nu2+delta_to_zeta2,delta_def_h} below,
	and $\int h^{}_{\delta} \,\dd \rmN \leq \frac{25}{24}\nu_{3}(\rmN) + \frac{4}{27}$.
\end{lemma}
\begin{proof}
	Let $f \coloneqq \abs{\,\cdot\,}^{2+\delta}$ and
	\begin{align} \label{Eq--nu2+delta_to_zeta2,delta_def_h}
	 h(x) \,\coloneqq\, 	h_\delta(x)  \, \coloneqq \,
	  \left.\begin{cases} \frac{\delta(2+\delta)}{3} \abs{x}^3
	   + \frac{(1-\delta)(2+\delta)}{2} x^2 & \ \text{if} \  \abs{x} \leq 1, \\ \abs{x}^{2+\delta}
	    + \frac{\delta(1-\delta)}{6} & \ \text{if} \ \abs{x} > 1
	  \end{cases}\right\} \quad \text{for $x \in \R$} \,.
	\end{align}
    If $\delta=0$, then $h_\delta=(\cdot)^2$ and the stated inequalities are trivial.
    Hence we assume $\delta>0$ from now on.

	The function $h$ is symmetric, positive, and twice differentiable with
	\begin{align} \label{Eq--nu2+delta_to_zeta2,delta_def_h''}
		h''(x)  \, = \,
		\left.\begin{cases} 2\delta(2+\delta) \abs{x} + (1-\delta)(2+\delta) & \ \text{if} \ \abs{x} \leq 1 \\ (1+\delta)(2+\delta) \abs{x}^{\delta} & \ \text{if} \ \abs{x} \ge 1 \end{cases}\right\} \quad \text{for $x \in \R$} \,.
	\end{align}
	Now let $x, y \in \R$. We intend to show that $\frac{1}{(1+\delta)(2+\delta)} h \in \cF_{2,\delta}$ by showing
	\begin{align} \label{Eq--nu2+delta_to_zeta2,delta_cond_to_sec_derivative}
		\abs{ h''(x) - h''(y) }  \, \leq \,  (1+\delta)(2+\delta) \left( \abs{x - y}^{\delta} \wedge \abs{x - y} \right) \,.
	\end{align}
	Since $h''(z) = h''(-z)$ for $z \in \R$, we may assume w.l.o.g.~$x, y \in \mathopen[0, \infty\mathclose[$. 
	From~\eqref{Eq--nu2+delta_to_zeta2,delta_def_h''} we see that $h''$ is continuous everywhere
	and differentiable on $\R\setminus\{-1,0,1\}$, with
	\begin{align}                            \label{Eq:h'''}
		h'''(z) \, = \, \sgn(z)\begin{cases} 2\delta(2+\delta) & \ \text{if} \ \abs{z} \in \mathopen] 0, 1\mathclose[\,,\\ \delta(1+\delta)(2+\delta)\abs{z}^{\delta - 1} & \ \text{if} \ \abs{z} \in \mathopen] 1, \infty \mathclose[   
		\end{cases} 
	\end{align}
	with $h'''(1-)= 2\delta(2+\delta) \ge \delta(1+\delta)(2+\delta) = h'''(1+)$.
	Hence we have
	\begin{align} \label{Eq--nu2+delta_to_zeta2,delta_|x-y|}
		\left| h''(x) - h''(y) \right|  \, \leq \,  \sup_{z \in \R\setminus\{-1, 0, 1\}} \left| h'''(z) \right|  \abs{x - y}  \, \leq \,  2\delta(2+\delta) \abs{x - y} \,.
	\end{align}
	On $[0,\infty[$ the function $h''$ is increasing and concave,
	by its continuity and by~\eqref{Eq:h'''}, 	and hence we get for $0\le x\le y<\infty$
	\begin{align*}
		\left| h''(y) - h''(x) \right|  \, = \,   h''(y ) - h''(x)
		 \, \leq \,   h''(y - x) - h''(0)   \, \le \, h''\big(|x - y|\big) \,,
	\end{align*}
	which yields 
	\begin{align} \label{Eq--nu2+delta_to_zeta2,delta_|x-y|>1}
		\left| h''(x) - h''(y) \right|  \, \leq \,  (1+\delta)(2+\delta) \abs{x-y}^{\delta} \quad \text{if $\abs{x-y} > 1$} \,. 
	\end{align}
	Combining~\eqref{Eq--nu2+delta_to_zeta2,delta_|x-y|} and~\eqref{Eq--nu2+delta_to_zeta2,delta_|x-y|>1} with the estimate $(1+\delta)(2+\delta)  \, \geq \,  2\delta(2+\delta)$
	then yields~\eqref{Eq--nu2+delta_to_zeta2,delta_cond_to_sec_derivative}. 
	
	Since  $h(0)-f(0)=0$ and
	$h'(x)-f'(x)= (2+\delta)x(1-\delta+\delta x-x^\delta)\ge 0$ for $x\in[0,1]$,
	by the tangent upper bound at $x=1$ for the concave function $(\cdot)^\delta$,
	we see that  $h\ge f$ on~$\R$. Hence we get for $P \in \widetilde{\cP_{2+\delta}}\,$, using
	Lemma~\ref{Lem:cF_m,g-bounds_and_zeta_m,g_for_unbounded_f}\ref{part:zeta_m,g_for_unbounded_f}
	in the second step,
	\begin{align*}
		\nu_{2+\delta}(P)  \, \leq \,  \int h \,\dd(P - \rmN) + \int h \,\dd \rmN  \, \leq \,  (1+\delta)(2+\delta) \zeta_{2,\delta}( P - \rmN) + \int h \,\dd \rmN \,.
	\end{align*}
	Finally, we have for $x \in \mathopen[-1 , 1\mathclose]$
	\begin{align*}
		h(x)  \, &\leq \,  \left( \frac{\delta(2+\delta)}{3} + \frac{(1-\delta)(2+\delta)}{2} \right) \abs{x}^3 + \frac{(1-\delta)(2+\delta)}{2} \cdot \frac{4}{27}  \, \leq \,  \frac{25}{24} \abs{x}^3 + \frac{4}{27} \,,
	\end{align*}
	using in the first step $x^2 \leq \abs{x}^3 + \frac{4}{27}$, and in the second
	that $2\delta(2+\delta)+3(1-\delta)(2+\delta) = 6+\delta(1-\delta)
	 \leq \frac{25}{4}$ and $(1-\delta)(2+\delta)= 2-\delta-\delta^2\leq 2$.
	 Since also $h(x)< |x|^3+\frac{1}{24} < \frac{25}{24} \abs{x}^3 + \frac{4}{27}$
	 for $|x|>1$, we get
	 $\int h^{}_{\delta} \,\dd \rmN \leq \frac{25}{24} \nu_{3}(\rmN) + \frac{4}{27}\,$.
\end{proof}

\begin{proof}[Proof of the claim ``$\zeta_2^\flat \,\le\, \zeta_2$
on $\cM_{2,2}$'' in~\eqref{Eq:zeta_2,delta_vs_zeta_2+delta}]
                                        \label{page:prof_of_zeta_2^flat_le_zeta}
Let $M\in\cM_{2,2}$. If $f\in\cF^\infty_{2,(\cdot)^0}$,
that is, if $f$ is twice differentiable with $|f''(x)-f''(y)|\le 1$
for $x,y\in\R$, then $\widetilde{f}(x)\coloneqq f(x)-f''(0)\frac{x^2}{2}$
for $x\in\R$ defines a function $\widetilde{f}$
twice differentiable with $|\widetilde{f}''(x)|=|f''(x)-f''(0)|\le 1$
for $x\in\R$, hence once differentiable with
$|\widetilde{f}'(x)-\widetilde{f}'(y)| \le |x-y|$ for $x,y\in\R$,
that is, $\widetilde{f}\in\cF_{1,(\cdot)^1}$.
Hence then
\begin{equation*}
 \left|\int f\,\dd M(x)\right|
  \,=\, \left|\int \widetilde{f}\,\dd M(x)\right|
  \,\le\, \zeta_2(M)
\end{equation*}
by using in the first step $\mu_2(M)=0$, and in the second
Lemma~\ref{Lem:cF_m,g-bounds_and_zeta_m,g_for_unbounded_f}\ref{part:zeta_m,g_for_unbounded_f}
or equivalently \cite[p.~57, (1.68) with $r=2$]{Mattner2024}.
\end{proof}
\begin{proof}[Proof of inequality~{\rm\eqref{Eq:23/27}} and 
its sharpness] \label{Proof-23/27}
Let $g\in\cG$ be normalised. Then we have
$g(u)\ge u\wedge 1 = g^{}_1(u)$ for $u\in[0,\infty[\,$.
For $x\in[0,\infty[\,$, we have $x^3\ge x^2-\frac{4}{27}\,$, with equality iff $x=\frac{2}{3}\,$.
Hence, for $P\in\cP_2\,$,
\begin{align*}                           \label{Eq:23/27}
	\nu_{2,g}(\widetilde{P}) &\,\ge\, \nu_{2,g^{}_1}(\widetilde{P})
	 \,=\, \int|x|^3\wedge x^2\,\dd \widetilde{P}(x)
	 \,>\, \int(x^2-\tfrac{4}{27})\,\dd \widetilde{P}(x)
	 \,=\, \tfrac{23}{27}\,,
\end{align*}
with equality in the limit at the strict inequality if
$\widetilde{P} =\frac{1-\epsilon}{2}\big(\delta_{-\frac{2}{3}}+ \delta_{\frac{2}{3}}\big)
+\frac{\epsilon}{2}\big(\delta_{-a_\epsilon} + \delta_{a_\epsilon})$
with $a_\epsilon \coloneqq \frac{1}{3}\sqrt{ \frac{5}{\epsilon} +4 }\rightarrow \infty$ for
$\epsilon \rightarrow 0$.
\end{proof}

The following lemma justifies our claim, near the end of section~\ref{sec:Intro},
that a scaling inequality stated in
\cite[p.~746, third display from below]{Senatov1980} is wrong,
even if some universal constant factor is allowed on its right hand side,
and also if, in addition, the inequality is only considered for one law being
the standard normal and the other one being arbitrarily close to it.
The enorm $\norm{\,\cdot\,}$
on $\cM$ involved in making this latter condition precise
could be, for example,
$\norm{\,\cdot\,}^{}_\mathrm{K}$ or $\zeta_1$,
with the latter satisfying the boundedness condition
by~\eqref{Eq:zeta_1-distance_to_normality_bounded}.

For $b\in\mathopen]0,\infty\mathclose[$ and $\kappa\in[0,1]$,
\begin{equation*}
	g^{}_{b,\kappa}(u) \,\coloneqq\, u\wedge \big(b+\kappa(u-b)_+\big)
	\,=\, u\wedge\big((1-\kappa)b+\kappa u\big) \quad \text{ for }u\in[0,\infty[
\end{equation*}
defines a function $g^{}_{b,\kappa}\in\cG$, and then
\begin{align*}
	\cG_0 \,\coloneqq\, \big\{\,  g^{}_{b,\kappa}  \,:\, b\in[1,\infty[\,,\ \kappa \in \mathopen]0,1\mathclose]\,\big\}
	\, \subseteq\, \cG
\end{align*}
consists of normalised and unbounded functions only.
\begin{lemma}                             \label{Lem:Senatov_scaling_inequality_wrong}
	Let $\epsilon >0$, and let $\norm{\,\cdot\,}$ be an enorm on $\cM$ 
	bounded on $\big\{ P -\mathrm{N} : P\in\widetilde{\cP_2}\big\}$. Then
	\begin{align*}
		\sup\left\{ \frac{ \zeta_{2,g}\big(\big(P-\mathrm{N}\big)(\frac{\cdot}{a}) \big) }
		{ \frac{a^2}{g(\frac{1}{a})}\big( \zeta_{2,g}\vee\norm{\,\cdot\,}\big)(P-\mathrm{N}) }
		\right. : \,\
		& g\in\cG_0\,,\ a\in\mathopen]0,1\mathclose]\,,\ P\in\widetilde{\cP_2}  \\
		&  \,\text{\rm with }0 < \big( \zeta_{2,g}\vee
		\norm{\,\cdot\,}\big)(P-\mathrm{N}) \le \epsilon\
		\left.\vphantom{ \frac{\zeta_{2,g}\big(\big(P-\mathrm{N}\big)(\frac{\cdot}{a})\big)}{\frac{a^2}{g(\frac{1}{a})}\big(\zeta_{2,g}\vee\norm{\,\cdot\,}\big)(P-\mathrm{N})} }\right\} \, = \, \infty\,.
	\end{align*}
\end{lemma}
\begin{proof}
	1. Let $b\in[1,\infty[\,$, $\kappa\in[0,1]\,$, $g\coloneqq g^{}_{b,\kappa}\,$, and
	\begin{align}                              \label{Eq:Def_f_b_in_Gegenbeispiellemma}
		f(x)&\,\coloneqq\,  f^{}_{b,\kappa}(x) \\
		 & \nonumber\,\coloneqq\, \left\{\!\begin{array}{l}
			\frac{1}{6}|x|^3   \\
			(1-\kappa)\big(\frac{b^3}{6} +\frac{b^2}{2}\big(\,|x|\!-\!b\, \big) + \frac{b}{2}\big(\,|x|\!-\!b \,\big)^2\big)
			+ \frac{\kappa}{6}|x|^3
		\end{array}\right\}
		\text{ for } |x|  \,\left\{\!\!\begin{array}{l} \le \\  > \end{array}\!\!\right\}\, b\,.
	\end{align}
	Then $f(0)=f'(0)=f''(0)=0$ and $f''(x) = g\big(|x|\big)$ for $x\in\R$, and
	in the notation of
	Lemma~\ref{Lem:cF_m,g-bounds_and_zeta_m,g_for_unbounded_f}
	hence $f= g^{(-2)} \circ\abs{\,\cdot\,} \le (\,\cdot\,)^{(-2)}\circ\abs{\,\cdot\,} =\frac{1}{6}\abs{\,\cdot\,}^3$,
	and by its parts~\ref{part:When_primitive of_g_in_cF_m,g} and  \ref{part:zeta_le_nu}
	then $f\in\cF_{2,g}$ and
	$\zeta_{2,g}(M) \,\le\, \int f\,\dd |M|$ for $M\in\cM_{2,2}\,$.
	
	\smallskip
	2. For $\lambda\in\mathopen]0,1\mathclose]$ and $t\in[1,\infty[$\,, let
	\begin{align*}
		Q \,\coloneqq\,& Q_t\,\coloneqq\,  \big(1-\tfrac{1}{t^2}\big) \delta_0
		+ \tfrac{1}{2t^2}\big( \delta_{-t}+\delta_{t}\big)  \,, \\
		P \,\coloneqq\,& P_{\lambda,t}\,\coloneqq\, (1-\lambda)\mathrm{N} + \lambda Q_t\,.
	\end{align*}
	Then $Q,P\in\widetilde{\cP_2}\,$, $P\neq\mathrm{N}$, and
	\begin{align}              \label{Eq:P-N_vs_Q-N}
		P-\mathrm{N}\,=\,\lambda\cdot(Q-\mathrm{N})\,.
	\end{align}
	
	We now assume that $a,\kappa\in\mathopen]0,1\mathclose]$ and $t\in[\frac{1}{a^2},\infty[$\,, and we
	apply step~1 to $b\coloneqq at$ to get 
	\begin{align}                             \label{Eq:Gegenbeispiel_zeta_le}
		\zeta_{2,g^{}_{b,\kappa}}(Q-\mathrm{N}) &\,\le\, \int f^{}_{b,\kappa}\,\dd\abs{Q-\mathrm{N}} \\ \nonumber
		&\,=\,  \tfrac{1}{t^2}f^{}_{b,\kappa}(t) + \int f^{}_{b,\kappa}\,\dd\mathrm{N}
		\, \le \, \tfrac{1}{6}\Big( \big( (1-\kappa)3a + \kappa\big)t + \nu_3(\mathrm{N}) \Big)
	\end{align}
	by using for the first summand just the definition~\eqref{Eq:Def_f_b_in_Gegenbeispiellemma} and
	$t\ge at = b$ and hence
	$\tfrac{1}{t^2}f^{}_{b,\kappa}(t) - \frac{1}{6}\kappa t
	 =\frac{1-\kappa}{6}ta(a^2-3a+3)
	 \le \frac{1-\kappa}{6}3at$,
	and for the second $ f^{}_{b,\kappa} \le  \frac{1}{6}\abs{\,\cdot\,}^3$.
	On the other hand, using $f^{}_{b,\kappa}\in\cF_{2,g^{}_{b,\kappa}}$ and
	Lemma~\ref{Lem:cF_m,g-bounds_and_zeta_m,g_for_unbounded_f}\ref{part:zeta_m,g_for_unbounded_f},
	we get
	\begin{align}                             \label{Eq:Gegenbeispiel_zeta_ge}
		\zeta_{2,g^{}_{b,\kappa}}\big(\big(Q-\mathrm{N}\big)(\tfrac{\cdot}{a}) \big)
		&\, \ge \, \int f^{}_{b,\kappa}(ax)\,\dd(Q-\mathrm{N})(x) \\  \nonumber
		&\, =\,  \tfrac{1}{t^2}f^{}_{b,\kappa}(at) - \int f^{}_{b,\kappa}(ax)\,\dd\mathrm{N}(x)
		\, >\, \tfrac{1}{6}\left(a^3t - \nu_3(\mathrm{N})\right) \,.
	\end{align}
	
	We apply the above with $(a,\kappa)$ fixed, but with $t\rightarrow\infty$,
	hence $b=at\rightarrow \infty$ and hence $g^{}_{b,\kappa}(\frac{1}{a})= \frac{1}{a}$ eventually,
	and with $\lambda=\lambda(t)\rightarrow0$ chosen, using~\eqref{Eq:P-N_vs_Q-N}, in such a way that we always have
	$\big( \zeta_{2,g^{}_{b,\kappa}}\vee\,\norm{\,\cdot\,}\big)(P_{\lambda,t}-\mathrm{N})
	=\lambda\cdot\big( \zeta_{2,g^{}_{b,\kappa}}\vee\,\norm{\,\cdot\,}\big)(Q_{t}-\mathrm{N})
	\le \epsilon$.
	Then the pairs $(P_{\lambda,t}, g^{}_{b,\kappa})$ may be taken for lower-bounding the supremum
	in our claim, but the fractions to be considered 
	may, by~\eqref{Eq:P-N_vs_Q-N} again, as well be taken at the pairs $(Q_t, g^{}_{b,\kappa})$.
	The lower limit of the latter fractions is,
	by (\ref{Eq:Gegenbeispiel_zeta_le}, \ref{Eq:Gegenbeispiel_zeta_ge})
	and by the boundedness condition on $\norm{\,\cdot\,}$,
	seen to be $\ge$~$a^3/\Big(a^2a\big( (1-\kappa)3a + \kappa\big)\Big) = 1/\big( (1-\kappa)3´a + \kappa\big)$.
	Hence the claim follows by taking  $\kappa$ and $a$  arbitrarily small.
\end{proof}

\section{Errata for \cite{Mattner2024}}  \label{sec:Errata}
When
\cite[p.~70, Theorem 3.1]{Mattner2024} was said there
to be essentially Zolotarev's (1986, 1997)~\cite{Zolotarev1986, Zolotarev1997},
and was hence called ``Zolotarev's $\zeta_1\vee\zeta_3$ Theorem'',
it was overlooked that Zolotarev~(1997)~\cite[p.~ix]{Zolotarev1997} writes
that ``V.\,V.~Senatov \ldots~prepared'' it. So it should perhaps
rather be called ``Senatov's and Zolotarev's $\zeta_1\vee\zeta_3$ Theorem''.
We accordingly describe here its partial generalisation
Theorem~\ref{Thm--Z/S_zeta1_zeta2,delta} as being of Senatov-Zolotarev type.

A result of Senatov~(1998)~\cite[p.~162, Theorem 4.3.1]{Senatov1998}
was cited in \cite[p.~59, reference for (1.80)]{Mattner2024}
with a wrong page number, and it was overlooked
that it was already announced in Senatov~(1980)~\cite[p.~749, Remark~3]{Senatov1980}.

In \cite[p.~69, Question 2.2]{Mattner2024}, the restriction ``$n\ge2$''
should be replaced with ``$2\le n \le N-2$'', since the case of $n=N$ makes no sense at all,
and the cases of $n=N-1$ and $n=1$ are equivalent to each other.

In \cite[p.~85, (5.17)]{Mattner2024}, $\lambda_1(M_1)$ should be replaced
with $\abs{\lambda_1(M_1)}$.

In \cite[p.~94, display above (6.5)]{Mattner2024}, $\zeta_1\big(\widetilde{P^{\ast n}},\mathrm{N}\big)$
should be replaced with $\zeta_1\big(\widetilde{P^{\ast n}}-\mathrm{N}\big)$.

In \cite[p.~106, Example~12.2(b)]{Mattner2024}, ``$\text{L.H.S(3.9)}\sim\frac{1}{2}\text{R.H.S(3.9)}$''
should be replaced with
``$\text{L.H.S.(3.9)}\sim\text{R.H.S.(3.9)}\sim\frac{1}{2}\text{R.H.S.(3.10)}$''.

In \cite[p.~57]{Mattner2024}, the sets $\cF_r$ and their subsets should have been defined
to consist of $\R$-valued functions only.
For the occasional use of $|\int g\,\dd M|\le \zeta_{r}(M)$ for appropriate $\C$-valued $g$,
as in \cite[p.~90, step 5]{Mattner2024}, the following should have been inserted somewhere between
\cite[pp.~83--85, Lemma 5.1 and Corollary 5.2]{Mattner2024}:

In view of the following lemma, one can usually define
$\left\|M\right\|_\cF\coloneqq \sup_{f\in\cF}\left|\int f\,\dd M \right|$ for $M\in\cM$
with a set $\cF$ of $\R$-valued functions,
but still use $\left\|\,\cdot\,\right\|_\cF$
as if it were defined instead with some analogous set $\cG$ of $\C$-valued functions.

\begin{lemma}
Let $\cF$ and $\cG$ be sets of bounded and measurable $\C$-valued functions on $\R$.
\begin{parts}
\item If $\cF\subseteq\cG$, then $\left\|\,\cdot\,\right\|_\cF \le \left\|\,\cdot\,\right\|_\cG\,$.
\item If the implication
\begin{equation}
  g\in\cG\,,\ c\in\C\,,\ |c|=1 \ \Rightarrow \ \Re(cg) \in \cF
\end{equation}
holds, then $\left\|\,\cdot\,\right\|_\cF \ge \left\|\,\cdot\,\right\|_\cG\,$.
\end{parts}
\end{lemma}
\begin{proof} (a) is trivial.
 For (b) see the paragraph in the present paper after \eqref{Eq:scaling_identity_zeta_s}.
\end{proof}

\section*{Acknowledgement}
We thank a reviewer and the Editors for several helpful suggestions.


\providecommand{\bysame}{\leavevmode\hbox to3em{\hrulefill}\thinspace}
\providecommand{\MR}{\relax\ifhmode\unskip\space\fi MR }
\providecommand{\MRhref}[2]{%
	\href{http://www.ams.org/mathscinet-getitem?mr=#1}{#2}
}
\providecommand{\href}[2]{#2}

\section*{Notes on the references}

In the present paper, any page numbers in a citation
of a Russian source refer to its English translation. Links 
are provided here only if they lead to the full text for
free
as in case of~\cite{Esseen1945}, or are promised to be free if in a volume older than five years as in case of~\cite{Mattner2024}, 
or are read-only free as in case of~\cite{KorolevDorofeeva2017},
as experienced by us in July 2025.
The asterisk *  marks a reference
we have taken from secondary sources, without looking at the original.


\begin{thebibliography}{1}

	\bibitem{Esseen1945}
	C.-G.~Esseen, \emph{Fourier analysis of distribution functions. A mathematical study of the Laplace--Gaussian law.}
	Acta Math.~\textbf{77} (1945), no.~1, 1--125,
	\url{https://doi.org/10.1007/BF02392223}\,.

	\bibitem{Esseen1956}
	\bysame, \emph{A moment inequality with an application to the central limit theorem.}
	Skandinavisk Aktuarietidskrift~\textbf{39} (1956), 160--170.

	\bibitem{Goldstein2010}
    L.~Goldstein, \emph{Bounds on the constant in the mean central limit theorem.}
    Ann.~Probab.~\textbf{38} (2010), no.~4,
	1672--1689,
	\url{https://doi.org/10.1214/10-AOP527}\,.
	
	\bibitem{Jonas2024}
	L.~Jonas, \emph{An estimate for the rate of convergence in the central limit
		theorem using Zolotarev's distances for probability measures},
	unpublished master's thesis,
	Universit\"at Trier, 2024.
	
	\bibitem{Katz1963}
	M.~L.~Katz, \emph{Note on the Berry-Esseen theorem}, Ann. Math.	Statist.~\textbf{34} (1963), no.~3, 1107--1108,
	\url{https://doi.org/10.1214/aoms/1177704037}\,.
	
	\bibitem{KorolevDorofeeva2017}
	V.~Korolev and A.~Dorofeeva, \emph{Bounds of the accuracy of the
	normal approximation to the distributions of random sums under relaxed moment
	conditions}, Lith. Math. J.~\textbf{57} (2017), no.~1,
	38--58,
	\urlro{https://rdcu.be/eczr9}\,.
	
	\bibitem{Lindeberg1922}
	J.~W.~Lindeberg, \emph{Eine neue Herleitung des Exponentialgesetzes in der  Wahrscheinlichkeitsrechnung.}
	Math. Z.~\textbf{15} (1922), 211--225,
	\urlro{https://rdcu.be/dtanj}.
	
	\bibitem{Lyapunov1901}
	A.~Lyapunov, \emph{Nouvelle forme du th\'eor\`eme sur la limite de probabilit\'e},
	Zapiski Imperatorsko\u{\i} akadem\={\i}i nauk, po Fiziko-matematicheskomu otdielen\={\i}iu. 8e s\'erie = 
	M\'emoires de l'Acad\'emie imp\'eriale des sciences de St.-P\'etersbourg, Classe physico-math\'ematique. 8e s\'erie~\textbf{12}
	(1901), no.~5, 1--24 (French),
	\url{https://www.biodiversitylibrary.org/page/57713987}\,.
	
	\bibitem{Matskyavichyus1983}
	V.~K.~Matskyavichyus,
	\emph{A lower bound for the convergence rate in the central limit theorem},
	Teor.\ Veroyatn.\ Primen.~\textbf{28} (1983), no.~3, 565--569 (Russian),
	\url{https://www.mathnet.ru/eng/tvp2199}\,.
	English transl.~Theory Probab. Appl.~\textbf{28}, 596--601.
    Addendum in Teor. Veroyatn. Primen.~\textbf{29} (1984), no.~1, 198 (Russian),
    \url{https://www.mathnet.ru/eng/tvp1993}\,.
    English transl. Theory Probab. Appl.~\textbf{28} (1985), no.~1, 196--197.

	\bibitem{Mattner2024}
	L.~Mattner, \emph{A convolution inequality, yielding a sharper Berry–Esseen
		theorem for summands Zolotarev-close to normal}, Theory Probab. Math. Statist.~\textbf{111} (2024), 45--122,
		\urlfa{https://doi.org/10.1090/tpms/1217}\,. 

	\bibitem{MattnerShevtsova2019} L.~Mattner and I.~Shevtsova,
	 \emph{An optimal Berry–Esseen type theorem for integrals of smooth functions},
		ALEA, Lat. Am. J. Probab. Math. Stat.~\textbf{16} (2019), 487--530,
	 \url{http://alea.impa.br/articles/v16/16-19.pdf}.

	\bibitem{OsipovPetrov1967}
	 L.~V.~Osipov and V.~V.~Petrov,
	 \emph{On an estimate of the remainder term in the central limit theorem},
	 Teor. Veroyatn. Primen.~\textbf{12} (1967), no.~2,  281--286 (Russian),
	\url{https://www.mathnet.ru/eng/tvp708}\,.
	 English transl. Theory Probab. Appl.~\textbf{12}, no.~2, 281--286.
	
	\bibitem{Petrov1995}
	V.~V.~Petrov, \emph{Limit theorems of probability theory. Sequences of
		independent random variables},
		Oxford University Press, 1995.
	
	\bibitem{Senatov1980}
	V.~V.~Senatov, \emph{Several uniform estimates of the rate of convergence in the multidimensional central limit theorem.}
	Teor. Veroyatn. Primen.~\textbf{25} (1980), no.~4, 757--770 (Russian),
	\url{https://www.mathnet.ru/eng/tvp1230}\,.
	English transl. Theory Probab. Appl.~\textbf{25} (1981), no.~4, 745--759.
	
	\bibitem{Senatov1998}
	\bysame, \emph{Normal approximation: New results, methods and problems},
	VSP, Utrecht, The Netherlands,
	1998.
	
	\bibitem{Shevtsova2010}
	I.~G.~Shevtsova, \emph{On the asymptotically exact constants in the Berry–Esseen–Katz inequality}, Teor. Veroyatn. Primen.~\textbf{55} (2010), no.~2, 271--304 (Russian),
	\url{https://www.mathnet.ru/eng/tvp4201}\,.
	English transl.\ and revision
	Theory Probab. Appl.~\textbf{55}, no.~2 (2011), 225--252.

	\bibitem{Shevtsova2017}
	\bysame, \emph{On the absolute constants in Nagaev–Bikelis-type inequalities},
	Inequalities and extremal problems in probability and statistics. Selected topics
    (I.~Pinelis, ed.), Academic Press, 2017, pp.~47--102.

    \bibitem{Tyurin2010}
    I.~Tyurin, \emph{On the convergence rate in Lyapunov's theorem.}
    Teor. Veroyatn. Primen.~\textbf{55} (2010), no.~2, 250--270,
    \url{http://mi.mathnet.ru/eng/tvp4200} (Russian).
    English transl.~Theory Probab. Appl.~\textbf{55} (2011), no.~2, 253--270.

    \bibitem{Tyurin2012}
	\bysame, 
	\emph{Some optimal bounds in the central limit theorem using zero biasing},
    Statist. Probab. Letters~\textbf{82} (2012), no.~3, 514--518.

	\bibitem{Zolotarev1976a}
	V.~M.~Zolotarev, \emph{Approximation of distributions of sums of independent random
	variables with values in infinite-dimensional spaces},
	Teor. Veroyatn. Primen.~\textbf{21} (1976), no.~4, 741--758,
	\url{https://www.mathnet.ru/eng/tvp3420} (Russian).
	English transl.~Theory Probab.\ Appl.~\textbf{21} (1977), no.~4,	721--737.
	Errata in Teor. Veroyatn. Primen.~\textbf{22} (1977), no.~4, 901,
	\url{https://www.mathnet.ru/eng/tvp3642} (Russian).
	English transl. Theory Probab. Appl.~\textbf{22} (1978), no.~4, 881.

	\bibitem{Zolotarev1976b}
	\bysame, 
	\emph{Metric distances in spaces of random variables and their distributions.}
	Math. USSR Sbornik~\textbf{30} (1976), no.~3, 373--401,
	\url{https://www.mathnet.ru/eng/sm2908}.
	
	\bibitem{Zolotarev1986}
 	$^\ast$\bysame. 
 	\textit{Contemporary theory of summation of independent random variables.}
 	Nauka, Moscow (Russian), 1986.

	\bibitem{Zolotarev1997}
	 \bysame. 
 	\emph{Modern theory of summation of random variables}.
 	VSP, Utrecht, The Netherlands,
 	1997.
\end{thebibliography}
\end{document}